\documentclass[11pt]{article}
 \usepackage{bm,mathtools}
  \usepackage{amsfonts}
 \usepackage[a4paper, margin=2cm]{geometry}
  \usepackage{enumitem}

\newcommand{\diver}{\operatorname{{div}}}
\newcommand{\curl}{\operatorname*{{curl}}}

\usepackage{makecell}

\usepackage{tikz}
\usetikzlibrary{arrows,decorations.markings,arrows.meta}
\usetikzlibrary{babel}
\usepackage{pgfplots}
\pgfplotsset{compat=1.16}

\usepackage{amsmath}
\usepackage{amssymb}

\usepackage[normalem]{ulem}

\newcommand{\be}{\begin{equation}}
\newcommand{\ee}{\end{equation}}

\newcommand{\x}{{\bm{x}}}
\newcommand{\y}{{\bm{y}}}
\newcommand{\z}{{\bm{z}}}
\newcommand{\rr}{\bm{r}}
\newcommand{\nn}{\bm{n}}

 \makeatletter\@addtoreset{equation}{section}\makeatother

\newcommand{\END}{\hfill$\Box$}

{\refstepcounter{theorem}\noindent%
\textbf{Remark \thetheorem.}}%
{\hfill $\Box$ \\[0.75ex]}
\usepackage{fancyhdr}
\title{Nystr\"om discretizations of boundary integral equations for the solution of 2D elastic scattering problems}
\date{June 22, 2022}
\author{V\'{\i}ctor Dom\'{\i}nguez\thanks{Dep. {Estadística, Informática y Matemáticas}, Universidad P\'{u}blica de Navarra. Campus de Tudela 31500 - Tudela, Spain, e-mail: victor.dominguez@unavarra.es.}\and Catalin Turc\thanks{  Department of
Mathematical Sciences, New Jersey  Institute of Technology,
Univ. Heights. 323 Dr. M. L. King Jr. Blvd, Newark, NJ 07102, USA, e-mail: catalin.c.turc@njit.edu.}}
\allowdisplaybreaks

\begin{document}
\maketitle
\begin{abstract}
  We present three high-order Nystr\"om discretization strategies of various boundary integral equation formulations of the impenetrable time-harmonic Navier equations in two dimensions. One class of such formulations is based on the four classical {Boundary Integral Operators ({BIO}s)} associated with the Green's function of the Navier operator. We consider two types of Nystr\"om discretizations of these operators, one that relies on Kussmaul-Martensen logarithmic splittings~\cite{chapko2000numerical,dominguez2021boundary}, and the other on Alpert quadratures~\cite{alpert1999hybrid}. In addition, we consider an alternative formulation of Navier scattering problems based on Helmholtz decompositions of the elastic fields~\cite{dong2021highly}, which can be solved via a system of boundary integral equations that feature integral operators associated with the Helmholtz equation. Owing to the fact that some of the {BIO}s that are featured in those formulations are non-standard, we use Quadrature by Expansion (QBX) methods for their high order Nystr\"om discretization. Alternatively, we use Maue integration by parts techniques to recast those non-standard operators in terms of single and double layer Helmholtz {BIO}s whose Nystr\"om discretizations is amenable to the Kussmaul-Martensen methodology. We present a variety of numerical results concerning the high order accuracy that our Nystr\"om discretization elastic scattering solvers achieve for both smooth and Lipschitz boundaries. We also present extensive comparisons regarding the iterative behavior of solvers based on different integral equations in the high frequency regime. Finally, we illustrate how some of the Nystr\"om discretizations we considered can be incorporated seamlessly into the Convolution Quadrature (CQ) methodology to deliver high-order solutions of the time domain elastic scattering problems.
   
\end{abstract}

\section{Introduction}

Boundary integral equation based numerical solutions of scattering problems
%has
{have} certain inherent advantages over their volumetric counterparts, and thus {they have}
%it has
attracted significant attention in the literature in the past four decades.  The dimensional reduction and implicit enforcement of radiation conditions that solvers based on BIE enjoy, however, come with the challenge of resolving singular boundary integrals. In the arena of Nystr\"om discretizations a powerful methodology to deal with such singularities is the logarithmic splitting technique of Kussmaul-Martensen, which consists of making explicit the singular parts of various boundary operator kernels and subsequently dealing with them on a case by case basis. This technique has been successfully applied to first produce Nystr\"om discretizations of the Helmholtz {Boundary Integral Operators (BIOs)}~\cite{Kress,KressH} and subsequently of Navier BIOs~\cite{chapko2000numerical,dominguez2021boundary}. In the case of elastic scattering problem the fundamental solution of the Navier problem is significantly more complicated than its Helmholtz counterpart, which brings a host of new challenges in performing the singularity splitting technique. As documented in~\cite{dominguez2021boundary}, the Kussmaul-Martensen singularity splitting technique, although effective, becomes extremely cumbersome in the case of Navier BIOs, which may be a deterrent for its practical use. We explore in this paper alternative high-order Nystr\"om discretization strategies that,  owing to their being simpler to implement, may be more attractive  than the Kussmaul-Martensen approach.

One such alternative strategy is the use of Alpert quadratures~\cite{alpert1999hybrid} which resolve to high order logarithmic singularities. Basically, as shown in~\cite{dominguez2021boundary}, the differences between the three elastodynamics BIO which are not themselves weakly singular and their elastostatic counterparts are all weakly singular, that is their most singular component is logarithmic, and thus the Alpert quadrature is directly applicable to the discretization of such difference operators. The most important observation is that Alpert quadratures are applicable in a black box manner to the discretization of those difference operators without the need whatsoever to explicitly account for the logarithmic part of those kernels. The elastostatic BIOs (sometimes referred to as elasticity BIOs), on the other hand, have been studied in detail~\cite{hsiao2008boundary}, and their Nystr\"om discretization is amenable to integration by parts and periodic Hilbert transform techniques~\cite{KressH}. Based on the strategy outlined above we implement high-order Alpert quadrature Nystr\"om discretizations of elastodynamics BIOs for both smooth as well as Lipchitz boundaries, incorporating in the latter case sigmoid graded meshes~\cite{Kress} to resolve corner singularities. Furthermore, following the roadmap in~\cite{labarca2019convolution,petropoulos2021nystr}, we extend our frequency domain elastodynamics solvers to high-order solutions of the time domain elasticity scattering problems via the Convolution Quadrature (CQ) methodology.  Indeed, using Laplace transforms, the CQ methodology~\cite{banjai2009rapid,banjai2010multistep,banjai2011runge} reduces the solution of retarded potential formulations of wave equations to the solution of ensembles of Laplace domain elastodynamic problems. Alpert quadratures are seamlessly applicable to the discretization of the elastodynamic BIOs featuring fundamental solutions of the Navier equation with complex frequencies which are required in the CQ methods, and thus we derive high-order in time solutions of time domain elastic scattering equations.

We also consider an entirely different BIE strategy for the solution of elastodynamics scattering problems that relies on the Helmholtz decomposition of the elastic waves into compressional and shear waves.
%~\cite{dong2021highly,lai2019framework}.
The scalar functions corresponding to this Helmholtz decomposition, in turn, are radiative solutions of the Helmholtz equation with the pressure and respectively shear wave number, whose normal and tangential traces are coupled on the boundary of the scatterer. Seeking those aforementioned scalar functions in the form of Helmholtz layer potentials of corresponding wave numbers, the Navier scattering problems is reduced to a system of BIE which feature Helmholtz BIOs. This approach has been advocated in~\cite{lai2019framework,dong2021highly} for Dirichlet boundary conditions, and we extend it in this work to the case of Neumann boundary conditions. Using Helmholtz layer potential representations, we arrive at a system of BIE that feature certain Helmholtz BIOs which, in addition to being expressed in terms of Hadamard finite parts integrals as they involve the Hessian of the Helmholtz Green function, are non-standard in the sense that they do not result from the application of the Dirichlet and Neumann boundary traces to the usual single and double layer Helmholtz potentials. Furthermore, the derivation of the ensuing system of BIE requires the use of certain jump conditions~\cite{kolm2003quadruple} which are not typically encountered in this context. This procedure allows us to recast non-standard Helmholtz BIOs in terms of double and single layer Helmholtz BIOs whose Nystr\"om discretization is relatively straightforward within the Kussmaul-Martensen paradigm.

We find that the Quadrature by Expansion (QBX) {method} is a more straightforward strategy~\cite{MR3106484,klockner2013quadrature} for the  discretization of the Helmholtz decomposition elasticity BIEs. QBX {relies}
%rely
on smooth extensions of layer potentials evaluated in the exterior/interior PDE domains onto the boundary which are achieved in practice via Fourier-Bessel expansions connected to the addition theorem for Hankel functions. Thus, the application of boundary traces to layer potentials amounts to term by term differentiation of the Fourier-Bessel expansions in the QBX framework. An attractive feature of QBX methods is that the evaluation of the Fourier-Bessel expansion coefficients does not require resolution of kernel singularities. As such, we derive relatively straightforward QBX Nystr\"om discretizations of the Helmholtz decomposition BIE which converge with high order for both smooth and Lipchitz scatterers through the use of panel Chebyshev meshes and Crenshaw-Curtis quadratures for the evaluation of the QBX coefficients.

Arguably, the discretization of the Helmholtz decomposition BIE formulation of the Navier scattering problems is simpler than that of the BIE counterparts that use the Navier Green functions. In order to gain more insight into the properties of these two types of BIE formulations, we undergo a comparative study on their iterative behavior in the high frequency regime. In the case of Navier Green function based BIE formulations for the solution of scattering problems,  both combined field and regularized combined field formulations are available in the literature~\cite{chaillat2015approximate,chaillat2020analytical,dominguez2021boundary}. These formulations are proven to be robust for all frequencies, and they exhibit superior iterative behavior in the high frequency regime to the Helmholtz decomposition BIE we consider in this paper, even after the combined field approach is applied to the latter.  The analysis of the robustness of the combined field Helmholtz decomposition BIE is currently under investigation. 

 The paper is organized as follows: in Section~\ref{setting} we introduce the Navier equations that govern elastodynamics waves in two dimensions and their fundamental solution, and we review their associated boundary layer potentials and integral operators associated; in Section~\ref{BIEf} we review several Combined Field BIE formulations of Navier scattering problems which were discussed in detail in our previous contribution~\cite{dominguez2021boundary} as well as the Helmholtz decomposition BIE derived in~\cite{dong2021highly} for Dirichlet boundary conditions and extended in this work to the case of Neumann boundary conditions;
 %we review in Section~\ref{hoCQ} Runge-Kutta CQ methods for the solution of elastic wave equations following the presentation in~\cite{betcke2017overresolving};
 several versions of Nystr\"om discretizations of the elastodynamics BIOs as well as the Helmholtz BIOs that are featured in the Helmholtz decomposition BIE are presented in Section~\ref{Nd}; finally, we present in Section~\ref{NR} a variety of numerical results showcasing the high-order convergence achieved by the Nystr\"om discretizations of elastodynamics frequency  as well as the  iterative behavior of various BIE formulations in the high frequency regime {; we  conclude Section~\ref{NR} with an application of the discretizations discussed in the present work to the solution of the transient elastic wave equation by multistep and multistage Convolution Quadrature methods.}

 %in Section~\ref{OS} we derive and analyze Optimized Schwarz formulations of elastodynamics transmission problems; and finally in Section~\ref{Nd} we present two versions of Nyst\"om discretizations of the elastodynamics BIOs as well as a variety of numerical results that illustrate the high-order convergence of these methods as well as the iterative behavior of the various formulations considered in this paper.
 
 %collect results related to logarithmic singularity splittings of Hankel functions which we then use in  Sections~\ref{sl}-\ref{w} to perform singularity splittings for the kernels of the four leastdynamics BIOs of the Calderon calculus; in Section~\ref{ps} we construct the principal symbols of the four BIOs in the sense of periodic pseudodifferential operators; in Section~\ref{BIEf} we construct and analyze regularized BIE formulations for both penetrable and impenetrable scattering problems; in Section~\ref{OS} we derive and analyze Optimized Schwarz formulations of elastodynamics transmission problems; and finally in Section~\ref{Nd} we present two versions of Nyst\"om discretizations of the elastodynamics BIOs as well as a variety of numerical results that illustrate the high-order convergence of these methods as well as the iterative behavior of the various formulations considered in this paper.

\section{Boundary elements methods for Navier equations}\label{setting}

\subsection{Navier equations}
Let ${\bf u}{(x_1,x_2)}=(u_1(x_1,x_2),u_2(x_1,x_2)):\mathbb{R}^2\to   \mathbb{R}^2$ be a vector function. For a linear isotropic and homogeneous elastic medium with Lam\'e constants $\lambda$ and $\mu$ such that $\lambda>-2\mu$, the strain   and stress tensor are {given by}
\begin{eqnarray*}
 \bm{\epsilon}({\bf u})&:=&\frac{1}2 (\nabla {\bf u}+(\nabla {\bf u})^\top)
 =\begin{bmatrix}
  \partial_{x_1} u_1& \tfrac{1}2\left(\partial_{x_1} u_2+\partial_{x_2} u_1\right)     \\
   \tfrac{1}2\left(\partial_{x_1} u_2+\partial_{x_2} u_1\right) &                                                    
   \partial_{x_2} u_2
     \end{bmatrix}\\
 \bm{\sigma}({\bf u})&:=&2\mu \bm{\epsilon}({\bf u})+\lambda ({\diver {\bf u}}) I_2,\quad
 \end{eqnarray*}
{where    $\diver {\bf u} :=\partial_{x_1}u_1+\partial_ {x_2}u_2$  is the divergence operator and  $I_2$,  obviously, the identity matrix of order 2.  The time-harmonic elastic wave (Navier) equation is defined by}
\[
 \diver\bm{\sigma}({\bf u})+\omega^2{\bf u}=\mu\Delta {\bf u}+(\lambda+\mu)\nabla(\diver {\bf u})+\omega^2 {\bf u}=0
\]
where the frequency $\omega\in\mathbb{R}^+$ and the divergence operator  {is applied to $\bm{\sigma}({\bf u})$ row-wise}. Considering a bounded domain $\Omega$ in $\mathbb{R}^2$ whose boundary $\Gamma$ is a closed
Lipchitz curve, we are interested in solving the impenetrable elastic scattering problem in ${\Omega_+}$, the exterior of $\Omega$, that is look for solutions of the time-harmonic Navier equation
\begin{equation}\label{eq:NavD}
  \diver\bm{\sigma}({\bf u})+\omega^2{\bf u}=0\quad{\rm in}\ {\Omega_+}:=\mathbb{R}^2\setminus\Omega
\end{equation}
that satisfy the Kupradze radiation condition at infinity {(cf. \cite{AmKaLe:2009}, \cite[Ch. 2]{KupGeBa:1979}; see also \eqref{eq:kupradze})}. On the boundary $\Gamma$ the solution ${\bf u}$ of~\eqref{eq:NavD} satisfies either the Dirichlet boundary condition
\[
 {\gamma_\Gamma {\bf u} ={{\bf u}|_{\Gamma}=}{\bf f}} %\quad{\rm on}\ \Gamma
\]
or the Neumann boundary condition
\[
T{\bf u}={\bf g} %\quad{\rm on}\ \Gamma
\]
where {
${\bf f}, {\bf g}:\Gamma\to \mathbb{C}$ are sufficiently regular functions  and $T$ is
the associated normal stress tensor (or traction operator)
on $\Gamma$ given by
\[
 T{\bf u}:=\bm{\sigma}({\bf u})\nn=\lambda (\diver  {\bf u}) \nn+2\mu (\nn\cdot {\nabla}){\bf u}-\mu(\curl{\bf u}){\bm t}.
\]
Here,  $\curl {\bf u}= \partial_{x_1} u_2-\partial_{x_2} u_1 $ is the rotational or scalar curl  of ${\bf u}$,   $\bm{n} $ the unit outward normal derivative and  $\bm{t} = (-n_2,n_1)$ is the unit positive orientated tangent vector.}
%
%
% In addition, we will consider the transient elastic wave equation which amounts to seeking the elastic wave
% \begin{equation}\label{eq:td}
% \partial_t^2 {\bf u}(\x;t)=\diver \bm{\sigma}( {\bf u}(\x;t))\qquad{\rm in}\ {\Omega_+}\times(0,\infty)
% \end{equation}
% with both Dirichlet
% \[
% {\bf u}(\x;t)={\bf g}(\x;t),\qquad \x\in\Gamma,\ t>0
% \]
% {\color{red}\large vD: Should we remove Neumann conditions?}
% or Neumann boundary conditions
% \[
% T{\bf u}(\x;t)={\bf f}(\x;t),\qquad \x\in\Gamma,\ t>0
% \]
%  where ${\bf g}(\x;t)$ and ${\bf f}(\x;t)$ are both causal, that is they vanish for $t<0$. We will use boundary integral equation formulations to solve the scattering problem~\eqref{eq:NavD} and the time domain problem~\eqref{eq:td} via Convolution Quadrature (CQ) methods. To this end, we begin by reviewing the fundamental solution of Navier equation and the associated layer potentials and {BIO}s.
%
% We will use boundary integral equation formulations to solve the scattering problem~\eqref{eq:NavD}. To this end, we begin by reviewing the fundamental solution of Navier equation and the layer potentials associated with the Navier equation.
%
\subsection{Fundamental solution of Navier equation and the associated boundary integral operators}

For  $\x=(x_1,x_2), \y=(y_1,y_2) \in\mathbb{R}^2$, we will denote 
\[
 \rr:=\x-\y,\quad r:=|\rr|=|\x-\y|.
\]
The fundamental solution of the time-harmonic elastic wave is given by
\begin{equation}\label{eq:defPhi}
 \Phi(\x,\y):=\Phi(\rr)=\frac{1}{\mu}\phi_0(k_s r)I_2+\frac{1}{\omega^2}\nabla_{\x}\nabla_{\x}^\top(\phi_0(k_s r)-\phi_0(k_p r)),\quad
  \phi_0(z) := \frac{\rm i}4 H_0^{(1)}(z),
\end{equation}
with $H_0^{(1)}$ the  Hankel function of first kind and order $0$ so that $\phi_0(kz)$ is the fundamental solution of the Helmholtz equation $\Delta \phi_0+k^2\phi_0=0$ and
\begin{eqnarray}
 k_p^2&:=&\frac{\omega^2}{\lambda+2\mu},\quad k_s^2:=\frac{\omega^2}{\mu}\label{eq:kp:ks}
 \end{eqnarray}
 the (squared) pressure and shear wave-number.
Boundary integral equation formulations of the Navier equations rely on the Navier layer potentials and their associated {BIO}s which we will review in what follows.
 
For a given density (vector) function ${\bm{\lambda}}:\Gamma\to\mathbb{C}^2$, the Navier single layer potential is defined as
\begin{equation}\label{eq:singleP}
  ({\bm{\mathcal{S}}}{\bm{\lambda}})(\z):=\int_\Gamma \Phi(\z,\y){\bm{\lambda}}(\y)\,{\rm d}\y,\quad \z\in\mathbb{R}^2\setminus\Gamma.
\end{equation}
The single layer potential is continuous in $\mathbb{R}^2$, and thus  the single layer boundary integral operator can be defined as
\begin{equation}\label{eq:single}
  ({\bm V}{\bm{\lambda}})(\x):=\lim_{\varepsilon\to 0}({\bm{\mathcal{S}}}{\bm{\lambda}})( \x+\varepsilon {\nn}(\x))=
  \int_\Gamma \Phi(\x,\y){\bm{\lambda}}(\y)\,{\rm d}\y ,\quad \x\in\Gamma.
  \end{equation}

%
% The Navier single layer BIO is weakly singular. Indeed, assuming a $2\pi$ parametrization ${{\bf x}}(t):\mathbb{R}\to\Gamma$,  the kernel of the parameterized single layer operator satisfies
% \[
%  \Phi({{\bf x}}(s),{{\bf x}}(t))=-\frac{1}{2\pi\omega^2}(k_s^2+k_p^2)\log \sin^2\tfrac{s-t}2  I_2
% + A(t,\tau)\sin^2\tfrac{s-t}2 \log {\sin^2\frac{t-s}2}+B(t,\tau),
% \]
% where $A$ and $B$ are $2\times 2$ bivariate $2\pi$ periodic smooth functions if ${{\bf x}}(t)$ itself is smooth.

{For   $ \bm{g} :\Gamma\to\mathbb{C}^2$, the double layer potential is, on the other hand, given by }
\begin{equation}\label{eq:doubleP}
  ({\bm{\mathcal{D}}} {\bm{g}})(\z):=\int_\Gamma \left[ T_{\y} \Phi(\z,\y)\right]^\top\bm{g}(\y)\,{\rm d}\y,\quad \z\in\mathbb{R}^2\setminus\Gamma
\end{equation}
where $T_{\y}\Phi(\z,\y)$ is {the} normal stress tensor applied column-wise to $\Phi(\z,\y)$ with respect to the  $\y$ variable. The double layer potential ${\bm{\mathcal{D}}}$ undergoes a jump discontinuity across $\Gamma$ so that{
% \[
% \lim_{\varepsilon\to 0}({\bm{\mathcal{D}}}{\bm{g}})(\x+\varepsilon {\nn}(\x))-\lim_{\varepsilon\to 0}({\bm{\mathcal{D}}}{\bm{g}})(\x-\varepsilon {\nn}(\x))={\bm{g}}(\x),\quad \x\in\Gamma
% \]
% and
\[
\lim_{\varepsilon\to 0^+}({\bm{\mathcal{D}}}{\bm{g}})(\x\pm \varepsilon {\nn}(\x))
%+\lim_{\varepsilon\to 0}({\bm{\mathcal{D}}}{\bm{g}})(\x-\varepsilon {\nn}(\x))
=\pm \frac12 \bm{g} + ({\bm{K}}\bm{g})({\x}),\quad \x\in\Gamma
\]
if $\Gamma$ is sufficiently smooth around ${\bf x}$, }
where the Double Layer BIO is defined explicitly as
\[
({\bm K}\bm{g})({\x}):={\mathrm{p.v.}}\int_\Gamma  {\bm{K}}(\x,\y)\bm{g}(\y)\,{\rm d}\y,\quad {\bm{K}}(\x,\y):= \left[ T_{\y} \Phi(\rr)\right]^\top.
\]
Unlike the single layer BIO, the double layer BIO ${\bm K}$ is no longer weakly singular and the integral above exists {only in the sense of a Cauchy principal value, hence the ``{p.v.}'' notation used above}.

%For a given density (vector) function ${\bm \varphi}$ defined on $\Gamma$,
The application of the traction operator to the single layer potential
gives rise to jump discontinuities
\[
\lim_{\varepsilon\to 0}(T {\bm{\mathcal{S}}}{\bm \varphi})(\x\pm \varepsilon {\nn}(\x))=\mp {\bm \varphi}(\x)
+({ \bm{K}^\top} \bm{\varphi})({\x}),\quad \x\in\Gamma
\]
where the adjoint double layer operator is given by 
\[
({\bm K}^\top \bm{\varphi})({\x}):=\int_\Gamma {\bm{K}}^\top (\x,\y)\bm{\varphi}(\y)\,{\rm d}\y,\quad {\bm{K}}^\top (\x,\y):=  T_{\x} \Phi(\x,\y).
\]

Finally, applying the traction operator %$T_{\x}$
to the double layer potential
%${\bm{\mathcal{D}}}{\bm{g}}$
we obtain
\[
\lim_{\varepsilon\to 0}T_{\x}({\bm{\mathcal{D}}}{\bm{g}})(\x\pm \varepsilon {\nn}(\x))=({\bm W} \bm{g})({\x}),\quad \x\in\Gamma
\]
where the BIO $W$ is defined as
\[
 ({\bm W} \bm{g})({\x}):={\rm f.p.}\int_\Gamma  {\bm{W}} (\x,\y)\bm{g}(\y)\,{\rm d}\y,\quad {\bm{W}} (\x,\y):=  T_{\y} \big[T_{\x} \Phi(\x,\y)\big].
\]
The kernel $W(\x,\y)$ is strongly singular (that is, it behaves like $\mathcal{O}(|\x-\y|^{-2})$ as $\y\to\x$), and as such the integral in its definition must be interpreted {in a Hadamard finite part sense which is precisely what ``{f.p.}'' stands for}.

\subsubsection{Singularity subtraction and parameterized version for the Navier Boundary Integral Operators}

Let us consider the  BIOs for elasticity $\bm{V}_0,\bm{K}_0, \bm{K}_0^\top$ and $\bm{W}_0
$, operators defined in the same manner from the fundamental solution of the elasticity problem
\[
 \Phi_0(\x,\y)=\frac{\lambda+3\mu}{4\pi\mu(\lambda+2\mu)}\left(-\log{r}\ I_2+\frac{\lambda+\mu}{\lambda+3\mu}{\bm G}(\rr)\right), \quad
 {\bm G}(\rr) =\frac{1}{r^2}\rr \rr^\top,
\]
It can be shown that the difference between the corresponding operators $\bm{K}-\bm{K}_0$, $\bm{K}^\top-\bm{K}_0^\top$ and $\bm{W}-\bm{W}_0$  are logarithmic, and so  weakly singular, integral operators.

Indeed, if ${\bf x}:\mathbb{R}\to \Gamma$ is a smooth, regular,  $2\pi-$periodic, counterclockwise oriented parameterization of a smooth curve $\Gamma$,  {introducing} the parameterized version of the densities
\[
 \bm{\lambda}(t) =\bm{\lambda}({\bf x}(t))|{\bf x}'(t)|,\quad
 \bm{g}(t) =\bm{g}({\bf x}(t)),
\]
and the  {parameterized normal stress tensor defined by}
\begin{equation}\label{eq:parameterizednormalstresstensor}
 T {\bf u}(t) =
(T {\bf u})({\bf x}(t)) |{\bf x}'(t)|
\end{equation}
we have the corresponding parameterized  BIO:
\begin{alignat*}{6}
   ({\bm V}\bm{\lambda})(t)&:=
  \int_0^{2\pi}\Phi({\bf x}(t  ),{\bf x}(\tau)){\bm{\lambda}}(\tau) \,{\rm d}\tau,\qquad &
   ({\bm K}\bm{g})(t) &:=  {\rm p.v.}\,
  \int_0^{2\pi} \left[T_{\tau}\Phi({\bf x}(t ),{\bf x}(\tau))\right]^\top {\bm{g}}(\tau) \,{\rm d}\tau\\
  ({\bm K}^\top \bm{\lambda})(t) &:={\rm p.v.}\,
  \int_0^{2\pi}T_{t}\Phi({\bf x}(t  ),{\bf x}(\tau)){\bm{\lambda}}(\tau) \,{\rm d}\tau,\qquad &
   ({\bm W}\bm{g})(t) &:=  {\rm f.p.}\,
  \int_0^{2\pi} T_t \left[T_{\tau}\Phi({\bf x}(t ),{\bf x}(\tau))\right]^\top {\bm{g}}(\tau) \,{\rm d}\tau.
\end{alignat*}
We then have that the kernels of the operators, denoted with a slight abuse of notation (we are confident that the context will make it clear whenever one of these operators arises whether it is the BIOs of the curve or its parameterized version)
by the same symbols ${\bm V}$, ${\bm K}$, ${\bm K}^\top$ and ${\bm W}$ can be factorized as
\begin{equation}\label{eq:BIOsparameterized}
\begin{aligned}
 {\bm V}(t,\tau)&\ = {\bm V}_0(t,\tau)
+ {\bm A}(t,\tau)\sin^2{\frac{t-\tau}2}\log {\sin^2\frac{t-\tau}2}+{\bm B}(t,\tau) \\
 {\bm K}(t,\tau)&\ = {\bm K}_0(t,\tau)
+ {\bm C}(t,\tau)\sin {\frac{t-\tau}2}\log \sin^2{\frac{t-\tau}2}+{\bm D}(t,\tau),\quad\
 {\bm K}^\top (t,\tau)\ =  ({\bm K}  (\tau,t) )^\top,\\
 {\bm W}(t,\tau)&\ = {\bm W}_0(t,\tau)+{\bm E}(t,\tau)\log \sin^2{\frac{t-\tau}2} +{\bm F}(t,\tau)\\
\end{aligned}
\end{equation}
where
\begin{eqnarray*}
 {\bm V}_0(t,\tau)&:=&-\frac{\lambda+3 \mu}{\mu(\lambda+2 \mu)} \frac{1}{4 \pi} \log r I_{2}+\frac{\lambda+\mu}{\mu(\lambda+2 \mu)} \frac{1}{4 \pi}
{\bm G}(\bm{r})\\
 {\bm K}_0(t,\tau)&:=&\frac{1}{2 \pi r^{2}}({\bf x}'(t) \cdot \boldsymbol{r})\left(\frac{\mu}{\lambda+2 \mu} I_{2}+2 \frac{\lambda+\mu}{\lambda+2 \mu}
{\bm G}(\boldsymbol{r})\right) \\
 {\bm W}_0(t,\tau)&:=&-\frac{\mu(\lambda+\mu)}{\lambda+2 \mu} \frac{\partial^{2}}{\partial \tau \partial t} \frac{1}{\pi}\left(-\log r I_{2}+
{\bm G}(\boldsymbol{r})\right)
\end{eqnarray*}
and ${\bm r}$ and $r$ are given now  by
\[
{\bm  r}={\bm r}(t,\tau)= {\bf x}(t)-{\bf x}(\tau)  ,  \qquad  r= |{\bm r}|, \qquad
{\bm G}({\bm r}) =\frac{1}{|r|^2} {\bm r}{\bm r}^\top.
\]
Let us notice that $ {\bm V}_0$,  $ {\bm K}_0$, $ {\bm K}_0^\top$  and $\bm{W}_0$   turn out to be the kernels of  the corresponding BIO for elasticity operator, which with the convention followed so far will be denoted, also, by the same symbols.
Functions ${\bm A}, {\bm B}, {\bm C}$,  ${\bm D}$ and ${\bm G}$ are smooth periodic matrix functions.
Besides, we notice that the derivation with respect the parameters $t,\ \tau$ is nothing but the weighted tangential derivative. That is, if
\[
(\partial_s {g})(\x) := (\nabla v(\x))\cdot\bm{t}({\bm{x}}),\quad \x \in\Gamma,\quad \text{with $\gamma_\Gamma v =g$}
\]
we have
\[
g'(t)= |{\bf x}'(t)| (\partial_s g)({\bf x}(t)).
\]
We will make  {extensive use} of the tangential derivative operator in the next sections.

As an interesting byproduct, the principal symbol of the operators can be easily derived from these expressions. We refer the interested reader to \cite{dominguez2021boundary} for more on this topic. Observe also that strongly singular part of $\bm{W}$ can be be rewritten as
\[
 (\bm{W}_0 {\bm g})(t) = \frac{\mu(\lambda+\mu)}{\lambda+2 \mu}\frac{1}{\pi} \int_0^{2\pi}  \partial_\tau  \left(-\log r I_{2}+
{\bm G}(\boldsymbol{r})\right) \bm{g}'(\tau)\,{\rm d}\tau
\]
in the parameterized space, which is   simply  the parameterized version of the well-known Maue-type formula for elasticity:
\begin{equation}\label{eq:maue}
 (\bm{W}_0 {\bm g})(\x) = \frac{\mu(\lambda+\mu)}{\lambda+2 \mu}\frac{1}{\pi} \int_{\Gamma} \partial_{s_{\bm{x}}} \left(-\log r I_{2}+
{\bm G}(\boldsymbol{r})\right) \partial_{s_{\bm{y}}}\bm{g} (\y)\,{\rm d}\y.
\end{equation}

\section{Boundary integral formulations}\label{BIEf}

We present in what follows various strategies to derive BIE formulations of elastic scattering problems. Besides the classical combined field formulations CFIE we consider regularized formulations that rely on the use of approximations of the {Dirichlet-to-Neumann (DtN in what follows)} operators whose analysis was given in our previous contribution~\cite{dominguez2021boundary}. The design of the regularized formulations for elastic scattering problems follows the blueprint from the Helmholtz case~\cite{turc1,dominguez2016well}.

\subsection{Combined field integral equations}

Just like in the Helmholtz case~\cite{BrackhageWerner}, the classical approach~\cite{chaillat2008fast,chaillat2017fast} in the case of Dirichlet boundary conditions is to look for a scattered field in the form of a Combined Field representation
\[
{\bf u}(\x):=({\bm{\mathcal{D}}}{{\bm{g}}})(\x)-i\eta_D ({\bm{\mathcal{S}}}{{\bm{g}}})(\x),\quad \x\in\mathbb{R}^2\setminus\Omega
\]
where the coupling parameter $\eta_D\neq 0$, leading to the Combined Field Integral Equation (CFIE)
\begin{equation}\label{eq:CFIE_D}
  \frac{1}{2}{{\bm{g}}}+{\bm K}{{\bm{g}}}-i\eta_D {\bm V}{{\bm{g}}}={{\bf f}}.
  %\quad {\rm on}\ \Gamma.
\end{equation}
The question of selecting a value of the coupling parameter $\eta_D$ that leads to formulations with superior spectral properties (and thus faster convergence rates for iterative solver solutions) can be settled via  DtN   arguments~\cite{chaillat2017fast}. We begin with Somigliana's identities
\[
{\bf u} = {\bm{\mathcal{D}}} {\gamma_\Gamma {\bf u}} -  {\bm{\mathcal{S}}}
{{T}{\bf u}} ,\quad  \text{in }{\Omega_+}= \mathbb{R}^2\setminus\Omega
\]
which we rewrite considering as $\gamma_\Gamma{\bf u}$ as the primary unknown boundary density and incorporating the DtN operator $Y$ in the form
\[
{\bf u} ={\bm{\mathcal{D}}} {\gamma_\Gamma {\bf u}} - {\bm{\mathcal{S}}}(
{ Y\gamma_\Gamma {\bf u}}  ),\quad \text{in }{\Omega_+} .
\]
The main idea in constructing regularized formulations is to use easy to construct approximations $\mathcal{R}^D$ of the DtN operator $Y$ and look for combined field representations in the form
\begin{equation}\label{eq:CFIER_D_potential}
  {\bf u} := {\bm{\mathcal{D}}}{\bm{g}} - {\bm{\mathcal{S}}}(\mathcal{R}^D{\bm{g}} )
  %(\x),\quad \x\in\mathbb{R}^2\setminus\Omega
\end{equation}
leading to the Combined Field Regularized Integral Equation (CFIER)
\begin{equation}\label{eq:CFIER_D}
  \frac{1}{2}{\bm{g}}+{\bm K}{\bm{g}}-{\bm V}\mathcal{R}^D{\bm{g}}={\bf f}. %\quad {\rm on}\ \Gamma.
\end{equation}
Clearly, the better the regularizing operator $\mathcal{R}^D$ approximates the DtN operator $Y$, the closer the operator in the left hand side of the CFIER equation~\eqref{eq:CFIER_D} is to the identity operator, and, in conclusion, the ensuing CFIER formulations are more suitable to iterative solver solutions (e.g. GMRES).  In our previous effort~\cite{dominguez2021boundary} we proposed the regularizing operator $\mathcal{R}^D= PS_\kappa(Y)$ which is a Fourier multiplier operator whose symbol $\sigma_\kappa[Y](\xi)$ is defined by
  \begin{equation}\label{eq:PSY}
  \sigma_{\kappa}[Y](\xi)= -\frac{4\mu(\lambda+2\mu)}{\lambda+3\mu}(|\xi|^2-\kappa^2)^{1/2}\left(\frac{1}{2}I_2+\frac{i\mu}{2(\lambda + 2\mu)}\begin{bmatrix}& \rm{sign}(\xi)\\ -\rm{sign}(\xi) & \end{bmatrix}\right)
  \end{equation}
where $\Im{\kappa}>0$. However, the regularizing operator $\mathcal{R}^D= PS_\kappa(Y)$ is a pseudodifferential operator of order $1$ whose numerical evaluation   can be  consequently more involved. Using the high-frequency approximation $|\kappa|\to\infty$ in equation~\eqref{eq:PSY} we can construct a simple regularizing operator
\[
\mathcal{R}^D_1=-\frac{2\mu(\lambda+2\mu)}{\lambda+3\mu}i\kappa,
\]
which, incidental, can be interpreted as delivering a quasi-optimal choice for the coupling parameter $\eta_D$ in the CFIE formulation
\begin{equation}\label{eq:eta_D}
\eta_D^{\rm opt}=\frac{2\mu(\lambda+2\mu)}{\lambda+3\mu}k_s
\end{equation}
if we choose $\kappa=k_s$. We remark that a similar, easily implementable low-order approximation of the DtN operator was proposed in~\cite{chaillat2017fast} as a regularizing operator in the Dirichlet case.

In the case of Neumann boundary conditions, we can look for a scattered field in the form of a Combined Field representation akin to the Burton-Miller formulation in the Helmholtz case~\cite{BurtonMiller}
\[
{\bf u} :=- {\bm{\mathcal{S}}}{\bm{\varphi}}  + i\eta_N {\bm{\mathcal{D}}}{\bm{\varphi}}   ,\quad \text{in } \Omega_+
\]
where the coupling parameter $\eta_N\neq 0$, leading to the Combined Field Integral Equation (CFIE)
\begin{equation}\label{eq:CFIE_N}
  \frac{1}{2}{\bm{\varphi}}-{\bm K}^\top {\bm{\varphi}}+i\eta_N {\bm W}{\bm{\varphi}}={{\bf g}}.%\quad {\rm on}\ \Gamma.
\end{equation}
Here again we start by recasting the Somigliana's identities looking at ${{T}{\bf u}} $ as the primary unknown boundary density and making use of the the Neumann-to-Dirichlet (NtD) operator (which is the inverse of the DtN operator)
\[
{\bf u}= {\bm{\mathcal{D}}} (Y^{-1}{{T}{\bf u}}) -  {\bm{\mathcal{S}}} {{T}{\bf u}} ,\quad \text{in } \Omega_+.
\]
The construction of regularized formulations relies again on available approximations $\mathcal{R}^N$ of the NtD operator $Y^{-1}$ via looking for combined field representations in the form
\begin{equation}\label{eq:CFIER_N_potential}
  {\bf u} := {\bm{\mathcal{D}}}(\mathcal{R}^N \bm{\varphi}  )- {\bm{\mathcal{S}}} \bm{\varphi}      ,\quad \text{in } \Omega_+
\end{equation}
leading to the Combined Field Regularized Integral Equation (CFIER)
\begin{equation}\label{eq:CFIER_N}
  \frac{1}{2}{\bm{\varphi}}-{\bm K}^\top{\bm{\varphi}}+{\bm W}\mathcal{R}^N\bm{\varphi}=
  {\bf g}.
\end{equation}
Again here, the choice $\mathcal{R}^N= PS_\kappa(Y^{-1})$ as the Fourier multiplier operator whose symbol is $(\sigma_\kappa[Y](\xi))^{-1}$ ( the inverse must be understood in matrix sense per formula~\eqref{eq:PSY}), leading again to well posed CFIER formulations, at least in the case when the boundary $\Gamma$ is smooth~\cite{dominguez2021boundary}.

Also, using high-frequency approximations we can construct a simple regularizing operator
\[
\mathcal{R}^N_1=i\frac{\lambda+3\mu}{2\mu(\lambda+2\mu)}\kappa^{-1},
\]
which, delivers a quasi-optimal choice for the coupling parameter $\eta_N$ in the CFIE formulation
\begin{equation}\label{eq:eta_N}
\eta_N^{\rm opt}=\frac{\lambda+3\mu}{2\mu(\lambda+2\mu)}k_s^{-1}.
\end{equation}
We remark that similar low-order approximations of NtD operators have been used in~\cite{chaillat2020analytical} to construct CFIE formulations with superior spectral properties.

Similarly, we can construct \emph{direct} regularized formulations in the case of Neumann boundary conditions following the ideas in~\cite{turc_corner_N}. Assuming that a smooth incident field ${\bf u}^{\rm inc}$ (which is a solution of the Navier equation in the whole $\mathbb{R}^2$) impinges on the obstacle $\Omega$, we will derive these BIEs in terms of unknown boundary quantity if ${\bf u}^{\rm tot}=({\bf u}+{\bf u}^{\rm inc})|_\Gamma$. We obtain from Somigliana's identities by taking into account the fact that ${{T}}{\bf u}^{\rm tot}=0$ on $\Gamma$
\begin{equation}\label{eq:S1}
  {\bf u} =  {\bm{\mathcal{D}}}{(\gamma_\Gamma{\bf u}^{\rm tot})}  ,\quad \text{in }\Omega_+.
\end{equation}
Applying the Dirichlet trace on $\Gamma$ to formula~\eqref{eq:S1} we obtain
\begin{equation}\label{eq:S11}
  \frac{1}{2}{(\gamma_\Gamma{\bf u}^{\rm tot})} - {\bm K}{(\gamma_\Gamma{\bf u}^{\rm tot})} ={\gamma_\Gamma}{\bf u}^{\rm inc}  ,
  %\quad \x\in\Gamma,
\end{equation}
while applying the traction operator to formula~\eqref{eq:S1} we get
\begin{equation}\label{eq:S12}
   {\bm W}{(\gamma_\Gamma{\bf u}^{\rm tot})} =-{{T}}{\bf u}^{\rm inc} .
\end{equation}
We combine BIE~\eqref{eq:S11} and a preconditioned (on the left) version of the BIE~\eqref{eq:S12} to arrive at the DCFIER
\begin{equation}\label{eq:S22}
  \frac{1}{2}%{\bf u}^{\rm tot}
  { \gamma_\Gamma{\bf u}^{\rm tot}}
  - {\bm K}
  {(\gamma_\Gamma{\bf u}^{\rm tot})}  +  \mathcal{R}^N{\bm W}
  {(\gamma_\Gamma{\bf u}^{\rm tot})}  =
  {\gamma_\Gamma{\bf u}^{\rm tot}} -\mathcal{R}^N({{T}}{\bf u}^{\rm inc}) .
\end{equation}
We note that the operators on the left hand side of the DCFIER formulation is the real $L^2(\Gamma)\times L^2(\Gamma)$ adjoint of the operator in the CFIER formulation, a situation that is, similar to that in the Helmholtz case~\cite{turc_corner_N}. We will use direct formulations in the case when $\Omega$ is a Lipschitz domain in order to take advantage of the increased regularity of ${\gamma_\Gamma}{\bf u}^{\rm tot}$.

\subsubsection{Open arcs}

In the case when $\Gamma$ is an open arc in $\mathbb{R}^2$, the combined field methodology is no longer available. Instead, first kind formulations can be derived from Somigliana's identities~\cite{chapko2000numerical}. Assuming again smooth incident fields ${\bf u}^{\rm inc}$, the scattering problem off of an arc $\Gamma$ is solved in the case of Dirichlet boundary conditions via the following BIE of the first kind
\begin{equation}\label{eq:arcS}
  {\bm V}[{{T}}{\bf u}^{\rm tot}]= \gamma_\Gamma {\bf u}^{\rm inc}
\end{equation}
while in the Neumann case via the BIE~\eqref{eq:S12}. The first kind integral equations of elastodynamic scattering from arcs can be preconditioned using the approximations of DtN operators introduced above. Thus, we will also consider the preconditioned BIE
\begin{equation}\label{eq:ParcS}
  (PS_\kappa(Y){\bm V})[T{\bf u}^{\rm tot}]= PS_\kappa(Y){\bf u}^{\rm inc}\quad {\rm on}\ \Gamma
\end{equation}
and
\begin{equation}\label{eq:ParcW}
  (PS_\kappa(Y^{-1}){\bm W})[{\bf u}^{\rm tot}]=-PS_\kappa(Y^{-1})[T{\bf u}^{\rm inc}]\quad {\rm on}\ \Gamma.
\end{equation}
We remark that a different strategy based on Calder\'on preconditioners have been proposed in~\cite{bruno2019weighted,bruno2020regularized} in order to produce formulations of elastodynamics scattering problems from arcs that have more suitable spectral properties to iterative solvers.

\subsection{Helmholtz decomposition formulations}\label{subsection:HmDeFor}

Another possibility to construct BIE formulations of Navier problems is via Helmholtz decompositions of the fields ${\bf u}$. Indeed, defining
\begin{equation}\label{eq:Hdecomp1}
{\bf u}_p:=-\frac{1}{k_p^2}\nabla \diver{\bf u}\qquad {\bf u}_s:=
%\frac{1}{k_s^2}\left[\nabla\curl}\ {\bf u}\right]^\perp=
{\frac{1}{k_s^2}\overrightarrow{\curl} \curl\, {\bf u}}
\end{equation}
with  $
   \overrightarrow{\curl}\,\varphi=[\nabla{\varphi}]^\perp=[\partial_{x_2}\varphi\ -\partial_{x_1}\varphi]^\top
$ the vector curl operator, we have that ${\bf u}={\bf u}_p+{\bf u}_s$. Hence, we can look for the fields ${\bf u}$ in the form
\begin{equation}\label{eq:Hdecomp2}
{\bf u}=\nabla\varphi_p+\overrightarrow{\curl} \ {\varphi_s}
\end{equation}
where the scalar functions $\varphi_p$ and $\varphi_s$ are radiative solutions of scalar Helmholtz equations in $\Omega_+$ with {wave-numbers} $k_p$ and $k_s$ respectively. In the case of Dirichlet boundary conditions, it is {simply} to see that $\varphi_p$ and $\varphi_s$ satisfy the following coupled boundary conditions
\begin{eqnarray}\label{eq:DH}
\partial_n\varphi_p+\partial_s\varphi_s&=&-{\gamma_\Gamma}{\bf u}^{\rm inc}\cdot \nn
%\qquad {\rm on}\ \Gamma
\nonumber\\
\partial_s\varphi_p-\partial_n\varphi_s&=&-{\gamma_\Gamma}{\bf u}^{\rm inc}\cdot \bm{t}.
\end{eqnarray}
Straightforward calculations {yield}
\begin{eqnarray*}
\bm{\sigma}(\nabla \varphi_p)&=&2\mu{\mathrm{H}}\varphi_p-\lambda k_p^2\varphi_p {\bm I}\\
\bm{\sigma}(\overrightarrow\curl\  \varphi_s)&=&2\mu{\mathrm{H}}\varphi_s\begin{bmatrix}&-1\\1 & \end{bmatrix}+\mu k_s^2\varphi_s \begin{bmatrix}&-1\\1 & \end{bmatrix}
\end{eqnarray*}
where ${\mathrm{H}}\varphi$ denotes the Hessian of the scalar function $\varphi$. {Using these calculations} in the case of Neumann boundary conditions we get in turn that $\varphi_p$ and $\varphi_s$ are coupled via the following boundary conditions
\begin{equation}
 \left\{
\begin{array}{rcl}
2\mu\ \nn^\top{\mathrm{H}}\varphi_p\nn-\lambda k_p^2\varphi_p+2\mu\ \nn^\top{\mathrm{H}}\varphi_s\bm{t}&=&-{T}{\bf u}^{\rm inc}\cdot \nn
%\qquad {\rm on}\ \Gamma
\\
2\mu\ \bm{t}^\top{\mathrm{H}}\varphi_p\nn-2\mu\ \nn^\top{\mathrm{H}}\varphi_s\nn-\mu k_s^2\varphi_s&=&-{T}{\bf u}^{\rm inc}\cdot \bm{t}. %\qquad {\rm on}\ \Gamma.
\end{array}
\right.
 \label{eq:NH}
\end{equation}
In both cases we look for $\varphi_p$ and $\varphi_s$ in the form of Helmholtz single layer potentials  corresponding to {wave-numbers} $k_p$ and respectively $k_s$. That is,  we look for  unknown functional densities ${g_p}$ and ${g_s}$ defined on $\Gamma$ such that
\begin{equation}\label{eq:SLrepr}
\varphi_p ={{\cal S}_{k_p}} {g_p}  \qquad
\varphi_s ={{\cal S}_{k_s}} {g_s}
\qquad{\rm in}\ \Omega_+,
\end{equation}
where {($\phi_0$ is the fundamental solution of the Helmholtz equation cf. \eqref{eq:defPhi})}
\[
 {{\cal S}_{k}g(\x):= \int_\Gamma \phi_0(k|\x-\y|)\,g(\y)\,{\rm d}y,\quad \x\in\Omega_+ }
\]
is the Single Layer BIO for Helmholtz equation with wave-number $k$ and. In the case of Dirichlet boundary conditions, the system of equations~\eqref{eq:DH} is equivalent to the following system of BIE for the boundary densities ${g_p}$ and ${g_s}$
\begin{equation}\label{eq:DHBIE}
{\bm{\mathcal{A}}_{\rm DH}^{\rm SL}}\begin{bmatrix}{g_p}\\  {g_s} \end{bmatrix}=-\begin{bmatrix}-{\bf u}^{\rm inc}\cdot \nn\\ -{\bf u}^{\rm inc}\cdot \bm{t}\end{bmatrix} \quad {\bm{\mathcal{A}}_{\rm DH}^{\rm SL}}=\begin{bmatrix} -\frac{1}{2}I+K_{k_p}^\top & \partial_s V_{k_s}\\ \partial_s V_{k_p} & \frac{1}{2}I -K_{k_s}^\top \end{bmatrix}
\end{equation}
where $V_k$ and $K_k^\top$ are the Helmholtz single layer and respectively the adjoint double layer BIOs associated with the Green's function $\phi_0(k\cdot)$ {and $\partial_s$ denotes the
tangential derivative operator on $\Gamma$}. The system of equations~\eqref{eq:DHBIE} was shown to be uniquely solvable in the case when $k_p^2$ and $k_s^2$ are not eigenvalues of $-\Delta$ in the interior domain $\Omega$ with Dirichlet boundary conditions~\cite{dong2021highly}. However, the analysis of the invertibility of the operator ${\bm{\mathcal{A}}_{\rm DH}^{\rm SL}}$ is quite involved on account of the degeneracy of its principal symbol in the pseudodifferential sense. Indeed, the principal symbol of the operator ${\bm{\mathcal{A}}_{\rm DH}^{\rm SL}}$ is a nilpotent matrix operator (actually its square equals zero), and thus the integral formulation~\eqref{eq:DHBIE} is far from an optimal formulation with regards to iterative solvers. Of course, it is possible to look for $\varphi_p$ and $\varphi_s$ in the form of Helmholtz double layer potentials corresponding to {wave-numbers} $k_p$ and respectively $k_s$
\begin{equation}\label{eq:DLrepr}
\varphi_p:={{\cal D}_{k_p}} {g_p} \qquad \varphi_s:={{\cal D}_{k_s}} {g_s} \qquad{\rm in}\ \Omega_+,
\end{equation}
where
\[
 {{\cal D}_{k} g (\x) =\int_\Gamma \partial_{\bm{n}(\y)}{\phi_0}(k|\x-\y|)g(\y)\,{\rm d}\y,
 \quad {\x}\in\Omega_+}
\]
is the Double Layer BIO for the Helmholtz Equation with {wave-number} $k$,
${g_p}$ and ${g_s}$ are unknown functional densities defined on $\Gamma$. In the case of Dirichlet boundary conditions, the system of equations~\eqref{eq:DH} is equivalent to the following system of BIE for the boundary densities ${g_p}$ and ${g_s}$
\begin{equation}\label{eq:DDHBIE}
{\bm{\mathcal{A}}_{\rm DH}^{\rm DL}}\begin{bmatrix}{g_p}\\  {g_s} \end{bmatrix}=-\begin{bmatrix}-{\bf u}^{\rm inc}\cdot \nn\\ -{\bf u}^{\rm inc}\cdot \bm{t}\end{bmatrix} \quad {\bm{\mathcal{A}}_{\rm DH}^{\rm DL}}:=\begin{bmatrix} W_{{k_p}} & \frac{1}{2}\partial_s+k_s^2{\bm t}\cdot V_{{k_s}}[\nn ]-K_{{k_s}}^\top\partial_s\\ \frac{1}{2}\partial_s+ k_p^2{\bm t}\cdot V_{{k_p}}[\nn ]-K_{{k_p}}^\top\partial_s& -W_{{k_s}}\end{bmatrix}
\end{equation}
where we denoted by $W_k$ the hyper singular Helmholtz BIO associated with {Helmholtz equation with wave-number $k$}.
%the Green's function $\phi_0(k\cdot)$.
We also took into account the Maue type formula~\cite{Kress} describing the behavior of the tangential derivative $\partial_s$ on $\Gamma$ applied to the exterior double layer potential, which we recount next
\begin{equation}\label{eq:tang_trace}
{\bm t}\cdot {\gamma_\Gamma(\nabla_+{{\cal D}_k\varphi})}=\frac{1}{2}\partial_s\varphi+k^2{\bm t}\cdot V_k[\varphi\nn]-K_k^\top[\partial_s\varphi],\
\end{equation}
The notation $\nabla_+$ in equation~\eqref{eq:tang_trace} refers to the application of the gradient in the exterior domain $\Omega_+$. It is not our intention to analyze the operator ${\bm{\mathcal{A}}_{\rm DH}^{\rm DL}}$ in what follows. Rather, we use a combined field approach
\begin{equation}\label{eq:CFIErepr}
\varphi_p:={{\cal D}_{k_p}} g_p -ik_p {{\cal S}_{k_p}} g_p \qquad
\varphi_s:={{\cal D}_{k_s}} g_s -ik_s {{\cal S}_{k_s}} g_s \qquad{\rm in}\ \Omega_+.
\end{equation}
leading to the CFIE formulation
\begin{equation}\label{eq:DHCFIE}
{\bm{\mathcal{A}}_{\rm DH}^{\rm CFIE}}\begin{bmatrix}g_p\\  g_s\end{bmatrix}=-\begin{bmatrix} {\bf u}^{\rm inc}\cdot \nn\\  {\bf u}^{\rm inc}\cdot \bm{t}\end{bmatrix} \quad {\bm{\mathcal{A}}_{\rm DH}^{\rm CFIE}}={\bm{\mathcal{A}}_{\rm DH}^{\rm DL}}-\begin{bmatrix}ik_p & \\ & ik_s\end{bmatrix}{\bm{\mathcal{A}}_{\rm DH}^{\rm SL}}
\end{equation}
which is more suitable for iterative solutions of elastic scattering problems. We leave the analysis of the well possedness of the CFIE formulation~\eqref{eq:DHCFIE} for future work.

In order to derive the Neumann counterpart of the BIE system~\eqref{eq:DHBIE} we rely on the following trace relation:
\begin{equation}\label{eq:traceH}
% \lim_{\x^+\to\Gamma}[\nabla \nabla_\top \ {\mathcal{S}_k}g](\x^+)={{\bf H}_{k}}[g](\x)+\frac{1}{2}\varkappa(\x)(-{\bm I}+ 2\nn \nn^\top)g(\x)-\frac{1}{2}(\nn{\bm t}^\top+{\bm t}\nn^\top)\partial_s g(\x)
 \gamma_\Gamma [{\mathrm{H}^+}  {\mathcal{S}_k}g] ={{\bf H}_{k}}g +\frac{1}{2}\varkappa (-{\bm I}+ 2\nn \nn^\top)g -\frac{1}{2}(\nn{\bm t}^\top+{\bm t}\nn^\top)\partial_s g.
\end{equation}
In the identity above,
{$\mathrm{H}^+$ means that the Hessian matrix operator is applied in $\Omega_+$,  the exterior of $\Gamma$, (to $\mathcal{S}_kg$)}
$\varkappa(\x)$ denotes the signed curvature of $\Gamma$ at $\x\in\Gamma$, and the matrix {BIO} ${{\bf H}_{k}}[g](\x)$ (understood in the sense of Hadamard finite parts) is defined as
\[
{{\bf H}_{k}}[g](\x):= {\mathrm{f.p.}}\int_\Gamma \nabla_\x \nabla_\x^\top \phi_0(\x-\y) g(\y) {\rm d}\y,\qquad \x\in\Gamma.
\]
We note that the formula~\eqref{eq:traceH} also appears in~\cite{klockner2013quadrature} and its justification relies on results established in~\cite{kolm2003quadruple}. We obtain the following system of BIE 
\begin{eqnarray}\label{eq:NH:02}
\mathcal{A}^{\rm SL}_{\rm NH}\begin{bmatrix}g_p\\  g_s \end{bmatrix}&=& -\begin{bmatrix} {T}{\bf u}^{\rm inc}\cdot \nn\\  {T}{\bf u}^{\rm inc}\cdot \bm{t}\end{bmatrix} \nonumber\\
 \mathcal{A}^{\rm SL}_{\rm NH}&=&
 \begin{bmatrix} 2\mu\ \nn^\top{\bm H}_{k_p}\nn& -\mu\partial_s \\
 -\mu\partial_s &-2\mu\ \nn^\top{\bm H}_{k_s}\nn
 \end{bmatrix}
 +
 \begin{bmatrix}
 \mu\kappa I -\lambda k_p^2V_{k_p} & 2\mu\bm{t}^\top {\bm H}_{k_s}\nn \\
 2\mu\ \bm{t}^\top{\bm H}_{k_p}\nn &-\mu\varkappa I-\mu k_s^2V_{k_s}
 \end{bmatrix}.
\end{eqnarray}
Indeed, formulas~\eqref{eq:NH} rely on the following jump relation for the Hessian of the single layer potential on $\Gamma$ which can be viewed as Maue's type formulas Hadamard finite parts integral operators into alternative expressions that involve Cauchy Principal Value and weakly singular integral operators that recast
\begin{eqnarray}\label{eq:jumpH}
\lim_{\x^+\to\x}{\mathrm{H}} {\cal S}_k [\varphi](\x^+)&=&-k^2\int_\Gamma \phi_0(\x-\y)\nn(\y)\nn^\top(\y)\varphi(\y){{\rm d}\y}\nonumber\\
&+&\int_\Gamma \left[\nabla_\x \phi_0(\x-\y){\bm t}^\top(\y)-{\bm Q}\nabla_\x \phi_0(\x-\y)\nn^\top(\y)\right]\partial_s\varphi(\y){{\rm d}\y}\nonumber\\
&-&\int_\Gamma \left[\nabla_\x \phi_0(\x-\y)\nn^\top(\y)+{\bm Q}\nabla_\x \phi_0(\x-\y){\bm t}^\top(\y)\right]\varkappa(\y)\varphi(\y){{\rm d}\y}\nonumber\\
&+&\frac{1}{2}\left[-I+2\nn(\x)\nn^\top(\x)\right]\varkappa(\x)\varphi(\x)-\frac{1}{2}\left[{\bm t}(x)\nn^\top(\x)+\nn(x){\bm t}^\top(\x)\right]\partial_s\varphi(\x).\nonumber\\
\end{eqnarray}

 We do not intend to provide neither a full derivation of formulas~\eqref{eq:jumpH}, nor an analysis of the integral formulation~\eqref{eq:NH} in this paper. It suffices to mention that it is straightforward to see that the principal symbol of the operators $ \mathcal{A}_{\rm NH}$ is again defective, and thus the analysis of the formulation~\eqref{eq:NH} requires a pseudodifferential calculus beyond the principal symbol. Also, the single layer formulation~\eqref{eq:NH} is not uniquely solvable for all material parameters. This situation can be remedied through the use of combined field formulations. However, the double layer analogue of the jump conditions~\eqref{eq:traceH} are significantly more involved, and we will devote a separate effort to their derivation in which we will present a complete analysis of robust integral formulations of Helmholtz decomposition based reformulations of elastic scattering problems. For the sake of completeness we simply recount the following jump relations
\begin{eqnarray}\label{eq:jumpHDL}
\lim_{\x^+\to\x}{\mathrm{H}} {\mathcal{D}_k}[\varphi](\x^+)
&=&k^2\int_\Gamma \nabla_\x \phi_0(\x-\y)\nn^\top(\y)\varphi(\y){{\rm d}\y}-\frac{k^2}{2}\nn(\x)\nn^\top(\x)\varphi(\x)\nonumber\\
&+&k^2\int_\Gamma \phi_0(\x-\y)\nn(\y)\bm{t}^\top(\y)\partial_s\varphi(\y){{\rm d}\y}\nonumber\\
&+&\int_\Gamma \left[\nabla_\x \phi_0(\x-\y)\nn^\top(\y)+{\bm Q}\nabla_\x \phi_0(\x-\y)\bm{t}^\top(\y)\right]\partial^2_s\varphi(\y){{\rm d}\y}\nonumber\\
&+&\int_\Gamma \left[\nabla_\x \phi_0(\x-\y)\bm{t}^\top(\y)-{\bm Q}\nabla_\x \phi_0(\x-\y)\nn^\top(\y)\right]\varkappa(\y)\partial_s\varphi(\y){{\rm d}\y}\nonumber\\
&-&\frac{1}{2}\left[{\bm Q}^\top+2\nn(\x){\bm t}^\top(\x)\right]\varkappa(\x)\partial_s\varphi(\x)-\frac{1}{2}\left[-I+2\nn(\x)\nn^\top(\x)\right]\partial^2_s\varphi(\x)\nonumber\\
\end{eqnarray}
where ${\bm Q}=\begin{bmatrix} & I\\ -I & \end{bmatrix}$.

{We then demote by $\mathcal{A}^{\rm DL}_{\rm NH}$ the boundary integral matrix operator derived from \eqref{eq:jumpHDL} which, for the sake of brevity do not detail here. Hence, the corresponding CFIER formulation, counterpart of is given by
\begin{equation}\label{eq:NHCFIE}
{\bm{\mathcal{A}}_{\rm NH}^{\rm CFIE}}\begin{bmatrix}g_p\\  g_s\end{bmatrix}=-\begin{bmatrix} T{\bf u}^{\rm inc}\cdot \nn\\  T{\bf u}^{\rm inc}\cdot \bm{t}\end{bmatrix} \quad {\bm{\mathcal{A}}_{\rm NH}^{\rm CFIE}}={\bm{\mathcal{A}}_{\rm NH}^{\rm DL}}-\begin{bmatrix}ik_p & \\ & ik_s\end{bmatrix}{\bm{\mathcal{A}}_{\rm NH}^{\rm SL}}.
\end{equation}
}

\section{ Nystr\"om discretizations}\label{Nd}

\subsection{Kussmaul-Martensen based discretizations}
We presented in a previous effort a Nystr\"om discretization strategy based on singularity splitting technique
{ and a resolution of logarithmic singularities which via the Kusmaul-Martensen quadrature~\cite{kusmaul,martensen}}.
In a nutshell, this discretization strategy relies on global trigonometric interpolation with $2n$ nodes
\[
%{\tau}
t_j=\frac{j\pi}{n},\ j=0,1,\ldots,2n-1
\]
onto the space of trigonometric polynomials
\[
T_n=\left\{\varphi(t)=\sum_{m=0}^na_m\cos{mt}+\sum_{m=1}^{n-1}b_m\sin{mt}\ :\ a_m,b_m\in\mathbb{C}\right\}
\]
and the use  of the Kusmaul-Martensen quadrature method, a product quadrature method based on the identity
\[
 -\frac{1}{2\pi}\int_{0}^{2\pi}\log (\sin^2(\tau/2))\cos m\tau\,{\rm d}\tau =\frac{1}{|m|},\quad m\in\mathbb{Z}\setminus\{0\}
\]
to resolve the logarithmic singularity.
We will refer to this method as K-M Nystr\"om method and is applicable to the kernels of the BIOs ${\bm V}$, ${\bm K}-{\bm K}_0$, ${\bm K}^\top-{\bm K}^\top_0$ and ${\bm W}-{\bm W}_0$). {With some modifications can be adapted to cover  Hilbert transform
quadratures~\cite{KressH}, as those appearing in ${\bm K}_0, {\bm K}_0^\top$,  and half grid size shifted quadrature methods~\cite{Kress} for the discretization of the static counterparts ${\bm W}_0$.} Given that the Nystr\"om discretizations of the principal parts of the elastodynamic operators are available in the literature, the main difficulty of the overall collocation schemes resides in the logarithmic splitting of various kernels. In other words, a precise description of the kernels of the elastodynamic operator presented in \eqref{eq:BIOsparameterized} is required which can be found in \cite[Appendix]{dominguez2021boundary}. Alternatively, we {will} present in next subsection a different collocation, based on specialized quadrature rules for logarithmic singular functions,  that acts in a black-box manner in the case of weakly singular kernels and thus bypasses the need for complicated kernel splittings.
 
 The Kusmaul-Martensen quadratures can be applied also for the Nystr\"om discretization of the Helm\-holtz decomposition formulations~\eqref{eq:DH} and~\eqref{eq:NH}. While the details of these were provided in the literature~\cite{dong2021highly} in the case of the single layer formulation~\eqref{eq:DH} Dirichlet boundary conditions, the Neumann boundary conditions counterpart~\eqref{eq:NH} is also amenable to such discretizations via the Maue type formulas~\eqref{eq:jumpH} and respectively~\eqref{eq:jumpHDL}. The details concerning the application of the Kusmaul-Martensen quadratures to  the operators in the
 equations~\eqref{eq:jumpHDL} and~\eqref{eq:jumpHDL} that feature the weakly singular kernels {$\phi_0(r)$ and $\partial_{\nn(\x)} \phi_0(r)$ (recall $r=|\x -\y|$)} can be found for instance in the classical reference~\cite{KressH}, while the operators that feature the singular kernel {$\partial_{s_{\x}} \phi_0(r)$} can be dealt with using the half grid shifting technique, see for instance~\cite{dominguez2021boundary}. Furthermore, the derivatives of the functional densities featured in equations~\eqref{eq:jumpH} and~\eqref{eq:jumpHDL} can be performed using Fourier differentiation.

\subsection{Alpert quadrature}

The guiding principle in the splitting calculations presented in the previous sections is the fact that the differences between dynamic and static versions of the elasticity BIOs are all weakly singular. However, those calculations become increasingly involved and cumbersome for the double layer and hypersingular elastodynamic BIO. We explore in what follows an alternative strategy for evaluations of those weakly singular operators that does not require complex singularity splittings. This alternative relies on Alpert quadratures~\cite{alpert1999hybrid,hao2014high} which can be applied seamlessly in our context. Specifically, Alpert quadrature takes on the form  
\begin{align*}
\int_0^{2\pi} k(t_i, \tau)\sigma(\tau)d\tau\ \approx\ & h\sum_{p=0}^{2n-2a} k(t_i,t_i+ah+ph)\sigma(t_i+ah +ph) \\
    &+ h\sum_{p=1}^m w_pk(t_i,t_i+\chi_p h)\sigma(t_i+\chi_p h ) \\
    & + h\sum_{p=1}^m w_p k(t_i,t_i+2\pi-\chi_p h)\sigma(t_i+2\pi-\chi_p h ),\qquad {h := \frac{2\pi}n}
\end{align*}
where the kernel $k(t_i,\tau)$ has a logarithmic singularity at $\tau= t_i$. Assuming that the density function $\sigma$ is a regular enough $2\pi$-periodic density, the integer parameter $a$, the weights $w_p$, and the nodes $\chi_p$ can be selected in such a matter so that the ensuing Alpert quadratures achieve prescribed high order convergence. The endpoint correction nodes $\chi_p$ are typically not integers, and as such the Alpert quadratures require evaluation of the density $\sigma$ outside of the equispaced mesh $\{t_i\}$. This can be achieved by resorting to Lagrange interpolation of order $m+3$ that shifts the grid points around the non-grid points where the density function $\sigma$ needs to be evaluated~\cite{alpert1999hybrid,hao2014high}. Specifically, we apply Alpert quadrature rules to the parametrized versions of the BIOs ${\bm V},\ {\bm K}-{\bm K}_0,\ {\bm K}^\top-{\bm K}_0^\top$ and ${\bm W}-{\bm W}_0$ {(cf. \eqref{eq:BIOsparameterized})}
without resorting to singularity splitting of their weakly singular kernels. The static BIOs (whose kernels correspond to the $0$ subindex in the notations above), on the other hand, are evaluated via the trigonometric interpolation Nystr\"om methods described above.
%We use in the numerical examples the 10-th order Alpert method (i.e. the parameter $a=6$ and $m=10$) in the case of smooth boundary curves $\Gamma$, and the 3-rd order Alpert method (i.e. the parameter $a=2$ and $m=3$) in the case of piece-wise smooth boundary curves $\Gamma$.

\subsection{Domains with corners}

The extension of the Nystr\"om discretizations of the elastodynamics BIOs to Lipschitz domains can be performed via sigmoid transforms that accumulate with algebraic orders discretization points toward corners points~\cite{dominguez2016well}. We assume that the domain ${\Omega}$ has corners at $\x_1,\x_2,\ldots,\x_P$ and that $\Gamma\setminus\{\x_1,\x_2,\ldots,\x_P\}$ is piecewise smooth. We assume that the boundary curve has a {regular $2\pi$-periodic smooth} parametrization so that each of the curved segments $[\x_j,\x_{j+1}]$ is parametrized by
\[
{\bf x}^{{w}}(t):=(x_1(w(t)), x_2(w(t))),\qquad
t\in[T_j,T_{j+1}]
\]
(so that $\x_j={\bf x}(T_j)$)  where $0=T_1<T_2<\ldots <T_P<T_{P+1}=2\pi$. Here
\[
w:[T_j,T_{j+1}]\to[T_j,T_{j+1}],\ 1\leq j\leq P
\]
is the sigmoid transform of order $p$ introduced by Kress~\cite{KressCorner}.  The function $w$ is a smooth and increasing bijection on each of the intervals $[T_j,T_{j+1}]$ for $1\leq j\leq P$, with $w^{(k)}(T_j)=w^{(k)}(T_{j+1})=0$
for $1\leq k\leq p-1$ and all $1\leq j \leq P$.
%With the aid of the graded meshes just introduced, we define \emph{weighted} traction traces on $\Gamma$ defined as
With the aid of the graded meshes, we can introduce the parameterized traction operator {as defined in \eqref{eq:parameterizednormalstresstensor}}
\[
{(T^w{\bf u})(\tau):={|{\bf x}'(w(\tau))w'(\tau)|}\ (T {\bf u})({\bf x}^{{w}}(\tau))}
\]
as well as  %\emph{weighted} parametrized BIOs acting on them
{the (weighted) parameterized {adjoint} double layer BIOs}
\begin{eqnarray*}
  %{\bm V}^w[T^w{\bf u}](\tau)&:=&\int_\Gamma \Phi({\bf x}(\tau),{\bf x}(t))\  [T^w{\bf u}]({\bf x}(t))\,{\rm d}t\\
 {
 {\bm K}^{\top,w}[T^w{\bf u}](t) := \int_0^{2\pi} T^w_{t}[\Phi({\bf x}^{{w}}(t),{\bf x}^{{w}}(\tau))]\  [T^w{\bf u}]( \tau)\,{\rm d}\tau.
 }
\end{eqnarray*}
{Since $|w'(\tau)|$ has zeros of order $p-1$ at the corners,  the singularity of the  parameterized density and the kernel of the adjoint double layer BIO   at the corners is cancel out}.

With regards to the double layer operator, we remark that while the kernel $K(\x,\y)-K_0(\x,\y)$ continues to be weakly singular in the Lipschitz case, the kernel $K_0(\x,\y)$ itself now becomes singular. We use a trick similar to that employed for the evaluation of the Laplace double layer operator in the case of domains with corners, that is, we recast the evaluation of the elastostatic double layer operator as
\[
\frac{1}{2}\bm{g}(\x) + ( {\bm K}_0 \bm{g})(\x)=\int_\Gamma {{\bm K}_0} (\x,\y)](\bm{g}(\y)-\bm{g}(\x))\ {\rm d}\y.
  \]
The previous identity is a direct consequence of the integration by parts formula (3.3.20) in~\cite{hsiao2008boundary} and well known results for the Laplace double layer operator~\cite{Kress}. For the evaluation of the hypersingular operator ${\bm W}$, we note that the operator ${\bm W}-{\bm W}_0$ is still weakly singular in the case when $\Gamma$ is Lipschitz. According to the discussion in Section~\ref{BIEf}, we prefer the use of direct formulations of scattering and transmission elastodynamic problems in Lipschitz domains, whose unknowns are ${\gamma_\Gamma {\bf u}}$ and/or ${T^w{\bf u}}] $ and which incorporate weighted BIOs. The Alpert discretizations can be immediately extended to the case of Lipschitz domains through the aid sigmoid transforms, an insight which has been first proposed in~\cite{labarca2019convolution} in the case of Helmholtz BIO. It is indeed an advantage of Alpert quadratures the feature that they are readily amenable to any situation which involves weakly singular kernels.

{Furthermore}, the Kusmaul-Martensen quadratures can be {also} extended in principle via sigmoid transforms for the discretization of the Helmholtz decomposition formulations~\eqref{eq:DH} and~\eqref{eq:NH} in the case of domains with corners.

\subsection{QBX}

Given that the BIOs featured in the BIE systems~\eqref{eq:DHBIE} and~\eqref{eq:NH} result from the application of various traces on $\Gamma$ to Helmholtz single layer potentials, QBX methods \cite{MR3106484,klockner2013quadrature}  are particularly advantageous alternatives in delivering high-order Nystr\"om discretizations of such BIOs, especially for the nonstandard ones involving Hessians of single layer potentials. Indeed, for a given {wave-number} $k$ and functional density $\varphi$ on $\Gamma$, the main thrust of QBX methods is extending single layer potentials ${{\cal S}}_k[\varphi](\x^+),\ \x^+\in\Omega_+$ to $\overline{\Omega_+}$ using Fourier-Bessel series expansions. The application of boundary traces to single layer potentials, therefore, is a matter of term by term differentiation of the Fourier-Bessel series expansions in the QBX discretization paradigm, as we explain in what follows. We start with presenting the details of the QBX method for the evaluation of single layer BIO. QBX relies on expansion centers, that is for $\x\in\Gamma$ we define $\x^\pm=\x\pm\varepsilon(\x) \nn(\x),\ \varepsilon(\x)>0$, and on the addition theorem for Hankel functions
\begin{equation}\label{eq:add_thm}
  H_0^{(1)}(k|\x-{\y}|)=\sum_{\ell=-\infty}^{\infty}H_\ell^{(1)}(k|\x^+-{\y}|)e^{i\ell(\theta')^+}J_\ell(k|\x-\x^+|)e^{-i\ell\theta^+},\ {\y}\in\Gamma
\end{equation}
where $\theta^+$ and $(\theta')^+$ are the angular coordinates of $\x$ and respectively ${\y}$ in the polar coordinate system centered at $\x^+$. Considering the smooth extension of ${{\cal S}}_k[\varphi](\x^+)$ onto $\Gamma$ as $\varepsilon(\x)\to 0$, we obtain the following Fourier-Bessel series representation of the single layer BIO $V_k[\varphi](\x)$
\begin{equation}\label{eq:qbx}
  V_k[\varphi](\x)=\sum_{\ell=-\infty}^{\infty}\alpha_\ell^+{(\x)} J_\ell(k|\x-\x^+|)e^{-i\ell\theta^+},\quad \alpha_\ell^+{(\x)}:=\frac{i}{4}\int_\Gamma H_\ell^{(1)}(k|{\x^+-\y}|)e^{i\ell(\theta')^+}\varphi({\y})\,{\rm d}{\y}.
\end{equation}
Since for any function $\psi$ defined on a neighborhood of $\Gamma$ the tangential and normal derivatives of the expansions~\eqref{eq:qbx}  in polar coordinates centered at $\x^+$  are given by
\[
 \partial_s \psi({\x}) =\nabla \psi({\x}
 )\cdot{\bm t}(\bm{x})= -\frac{1}{|\bm{x}-\bm{x}^+|}\partial_\theta\psi(r,\theta),\quad
 \partial_n \psi ({\x}) =\nabla \psi({\x})\cdot{\bm n}(\bm{x})= \partial_r\psi(r,\theta),\quad
\]
(here $
\psi(r,\theta ) = \psi(\x^++ ( r \cos\theta,r\sin\theta)),$ with $r = |\x -\x^+|$) we  have
\begin{equation}\label{eq:qbxDS}
  \partial_sV_k[\varphi](\x)=i\sum_{\ell=-\infty}^{\infty}\ell\alpha_\ell^+ \frac{J_\ell(k|\x-\x^+|)}{|\x-\x^+|}e^{-i\ell\theta^+}
\end{equation}
and
\begin{equation}\label{eq:qbxKtop}
  -\frac{1}{2}\varphi(\x)+K^\top_k[\varphi](\x)=-k\sum_{\ell=-\infty}^{\infty}\alpha_\ell^+({\x}) \left(-J_{\ell+1}(k{|\x-\x^+|})+\frac{\ell}{k{|\x-\x^+|}}J_{\ell}(k{|\x-\x^+|})\right)e^{-i\ell\theta^+},
\end{equation}
that is, the Fourier-Bessel series representations of the BIOs featuring in equations~\eqref{eq:DHBIE}. Notice that in \eqref{eq:qbxKtop} we have used the well-known identity for the derivative of the Bessel functions
\[
 J'_{\ell}(x)=-J_{\ell+1}(x)+\frac{\ell}x J_\ell(x), \quad\ell\in\mathbb{Z}.
\]
%
% in polar coordinates
% taking tangential and normal derivatives of the expansions~\eqref{eq:qbx} in polar coordinates centered at $\x^+$ we have
% and taking into account the continuity properties of gradients of single layer potentials across $\Gamma$ we derive Fourier-Bessel series representations of the BIOs featuring in equations~\eqref{eq:DHBIE}, first
% \begin{equation}\label{eq:qbxDS}
%   \partial_sV_k[\varphi](\x)=i\sum_{\ell=-\infty}^{\infty}\ell\alpha_\ell^+ \frac{J_\ell(k|\x-\x^+|)}{|\x-\x^+|}e^{-i\ell\theta^+}
% \end{equation}
% and then
% \begin{equation}\label{eq:qbxKtop}
%   -\frac{1}{2}\varphi(\x)+K^\top_k[\varphi](\x)=-k\sum_{\ell=-\infty}^{\infty}\alpha_\ell^+ \left(-J_{\ell+1}(k{|\x-\x^+|})+\frac{\ell}{k{|\x-\x^+|}}J_{\ell}(k{|\x-\x^+|})\right)e^{-i\ell\theta^+}.
% \end{equation}
Similarly, applying the Hessian in polar coordinates on the expansions~\eqref{eq:qbx}, taking into account appropriate jump relations on $\Gamma$ as $\varepsilon(\x)\to0$ and using instead
\[
 J'_{\ell}(x)=J_{\ell-1}(x)-\frac{\ell}x J_\ell(x), \quad\ell\in\mathbb{Z}
\]
we obtain the following Fourier-Bessel series representations of the BIOs featuring in equations~\eqref{eq:NH}, that is
\begin{eqnarray}\label{eq:qbxH1}
  \frac{\kappa(\x)}{2}\varphi(\x)&+&\nn^\top(\x) {\bf H}_k[\varphi](\x)\nn(\x)=\nonumber\\
&& \hspace{-1cm} \sum_{\ell=-\infty}^{\infty}\alpha_\ell^+{(\x)} \left(\frac{kJ_{\ell+1}(k{|\x-\x^+|})}{{|\x-\x^+|}}+\frac{(\ell^2-\ell-k^2{|\x-\x^+|}^2)J_{\ell}(k{|\x-\x^+|})}{{|\x-\x^+|}^2}\right)e^{-i\ell\theta^+}
  \end{eqnarray}
and respectively
\begin{eqnarray}\label{eq:qbxH2}
  -\frac{1}{2}\partial_s\varphi(\x)&+&\nn^\top(\x) {\bf H}_k[\varphi](\x){\bm t}(\x)=\nonumber\\
  &&i\sum_{\ell=-\infty}^{\infty}\ell\alpha_\ell^+{(\x)} \left(\frac{kJ_{\ell+1}(k{|\x-\x^+|})}{{|\x-\x^+|}}+\frac{(1-\ell)J_{\ell}(k{|\x-\x^+|})}{{|\x-\x^+|}^2}\right)e^{-i\ell\theta^+}.
  \end{eqnarray}
 It is also possible to average limits from $\Omega_+$ and ${\Omega}$ to obtain alternative series expansions for traces on $\Gamma$ of single layer potentials, especially in the case when those undergo jump discontinuities across the boundary. Specifically, we denote by $\theta^-$ and $(\theta')^-$ the the angular coordinates of $\x$ and respectively $\x'$ in the polar coordinate system centered at $\x^-$, and we derive the alternative Fourier-Bessel series expansion of the single layer BIO via smooth extensions from the interior domain $\Omega^{-}$ onto $\Gamma$:
\begin{equation}\label{eq:qbxm}
  V_k[\varphi](\x)=\sum_{\ell=-\infty}^{\infty}\alpha_\ell^{-}{ (\x)} J_\ell(k|\x-\x^{-}|)e^{-i\ell\theta^{-}},\quad \alpha_\ell^-(\x):=\frac{i}{4}\int_\Gamma H_\ell^{(1)}(k|\x'-\x^-|)e^{i\ell(\theta')^-}\varphi(\x')\,{\rm d}\x'.
\end{equation}
Using the jump conditions of traces of single layer potentials, we derive the two-sided Fourier-Bessel series representations
\begin{eqnarray}\label{eq:qbxKtopm}
  %-\frac{1}{2}\varphi(\x)+
  K^\top_k[\varphi](\x)&=&-\frac{1}{2}\varphi(\x)-\frac{k}{2}\sum_{\ell=-\infty}^{\infty}\alpha_\ell^+{(\x)} \left(-J_{\ell+1}(k{|\x-\x^+|})+\frac{\ell}{k{|\x-\x^+|}}J_{\ell}(k{|\x-\x^+|})\right)e^{-i\ell\theta^+}\nonumber\\
  &+&\frac{k}{2}\sum_{\ell=-\infty}^{\infty}\alpha_\ell^{-}{(\x)} \left(-J_{\ell+1}(k{|\x-\x^-|})+\frac{\ell}{k{|\x-\x^-|}}J_{\ell}(k{|\x-\x^-|})\right)e^{-i\ell\theta^{-}},
\end{eqnarray}
{we easily derive
\begin{eqnarray}\label{eq:qbxH1m}
  %\frac{\kappa(\x)}{2}\varphi(\x)&+&
  \nn^\top(\x) {\bf H}_k[\varphi](\x)\nn(\x)&=&
  %\frac{\kappa(\x)}{2}\varphi(\x)
  \nonumber
  \\
&&\hspace{-2cm}+  \frac{1}{2}\sum_{\ell=-\infty}^{\infty}\alpha_\ell^+{(\x)} \left(\frac{kJ_{\ell+1}(k{|\x-\x^+|})}{{|\x-\x^+|}}+\frac{(\ell^2-\ell-k^2{|\x-\x^+|}^2)J_{\ell}(k{|\x-\x^+|})}{{|\x-\x^+|}^2}\right)e^{-i\ell\theta^+}\nonumber\\
&&\hspace{-2cm}+\frac{1}{2}  \sum_{\ell=-\infty}^{\infty}\alpha_\ell^{-}{(\x)} \left(\frac{kJ_{\ell+1}(k{|\x-\x^-|})}{{|\x-\x^-|}}+\frac{(\ell^2-\ell-k^2{|\x-\x^-|}^2)J_{\ell}(k{|\x-\x^-|})}{{|\x-\x^-|}^2}\right)e^{-i\ell\theta^{-}}\qquad
  \end{eqnarray}
and finally
\begin{eqnarray}\label{eq:qbxH2m}
%  -\frac{1}{2}\partial_s\varphi(\x)&+&
\nn^\top(\x) {\bf H}_k[\varphi](\x){\bm t}(\x)&=&%-\frac{1}{2}\partial_s\varphi(\x)
\nonumber\\
  &&\hspace{-2cm}+\frac{i}{2}\sum_{\ell=-\infty}^{\infty}\ell\alpha_\ell^+{(\x)} \left(\frac{kJ_{\ell+1}(k{|\x-\x^+|})}{{|\x-\x^+|}}+\frac{(1-\ell)J_{\ell}(k{|\x-\x^+|})}{{|\x-\x^+|}^2}\right)e^{-i\ell\theta^+}\nonumber\\
  &&\hspace{-2cm}+\frac{i}{2}\sum_{\ell=-\infty}^{\infty}\ell\alpha_\ell^{-}{(\x)} \left(\frac{kJ_{\ell+1}(k{|\x-\x^-|})}{{|\x-\x^-|}}+\frac{(1-\ell)J_{\ell}(k{|\x-\x^-|})}{{|\x-\x^-|}^2}\right)e^{-i\ell\theta^{-}}.
  \end{eqnarray}}
  
  In the case when Double layer potentials are used in the Helmholtz decomposition formulations, their related BIOs that feature in the systems of boundary conditions~\eqref{eq:DH} and~\eqref{eq:NH} can be evaluated using the same term by term differentiation of Fourier-Bessel series strategy presented above simply replacing the coefficients $\alpha_\ell^+$ by their double layer counterparts defined as
  \begin{eqnarray}\label{eq:qbx1}
   \alpha_\ell^{{\rm DL},+}{(\x)}\!\!&:=&\!\!\frac{i}{4}\int_\Gamma \left(H_{\ell+1}^{(1)}(k|\x^+-{\y}|)-\frac{\ell}{k|\x^+-{\y}|}H_{\ell}^{(1)}(k|\x^+-{\y}|)\right)\frac{(\x^+-{\y})\cdot \nn({\y})}{|\x^+-{\y}|}e^{i\ell(\theta')^+}\varphi({\y}) {\rm d}{\y}\nonumber\\
   &-&\frac{\ell}{4}\int_\Gamma H_{\ell}^{(1)}(k|\x^+-{\y}|)\frac{(-\sin{(\theta')^+},\cos{(\theta')^+})\cdot \nn({\y})}{|\x^+-{\y}|}e^{i\ell(\theta')^+}\varphi({\y}){\rm d}{\y}.
\end{eqnarray}
The Fourier-Bessel expansions  above constitute the basis of QBX Nystr\"om discretizations of BIOs in equations~\eqref{eq:DHBIE} and~\eqref{eq:NH}. The full discretizations of those BIOs is achieved by (a) selecting a  truncating parameter $p$ in the Fourier-Bessel series above,
and (b) projecting the functional densities $\varphi$ into appropriate discrete functional spaces and thus effecting corresponding collocation quadratures for the evaluation of the expansion coefficients $\alpha_\ell^\pm$ for $-p\leq \ell\leq p$. We note that the integrands in the definition of coefficients $\alpha_\ell^\pm$ do not exhibit kernel singularities. In the case of smooth boundaries $\Gamma$, we consider global trigonometric interpolation of the densities $\varphi$ using the $2n$ equispaced nodes $t_m=\frac{m\pi}{n}, 0\leq m\leq 2n-1$ and trapezoidal quadratures for the evaluation of Fourier-Bessel expansion coefficients $\alpha_\ell^\pm({\bf x}(t_m)),\ 0\leq m\leq 2n-1$.
{In order to achieve uniform errors using the quadratures~\eqref{eq:trap} for all indices $\ell$ (the integrands in equations~\eqref{eq:trap} get increasingly oscillatory and the near-singularities of the Hankel functions more stringent as the indices $\ell$ get larger), Fourier interpolation is used to oversample the density $\varphi$ on a finer uniform mesh of size $2n'= 2\beta n$ where $\beta\in\mathbb{Z},\ \beta>1$. We then define
${\tau}_{j} = \frac{j\pi}{2 n'}$
and the coefficients $\alpha_\ell^\pm$ are evaluated using the trapezoidal rule on the finer mesh:
\begin{equation}\label{eq:trap}
\alpha_\ell^\pm({\bf x}(t_m))\approx \alpha_{\ell,n'}^\pm({\bf x}(t_m)):= \frac{i \pi}{4n'}\sum_{j=0}^{2n'-1}H_\ell^{(1)}(k|{{\bf x}^\pm}(t_m)-{\bf x}(\tau_j) |)e^{i\ell (\theta')^\pm}\varphi({\bf x}(\tau_j))|\x'(\tau_j)|
\end{equation} 
where ${{\bf x}^\pm}(t_m)={\bf x}(t_m)\pm\varepsilon_m \nn({\bf x}(t_m))$ with $\varepsilon_m=\min(|{\bf x}(t_m)-{\bf x}(t_{m-1})|,|{\bf x}(t_m)-{\bf x}(t_{m+1})|)$ (here we use cyclical indexing so that $t_{-1}:=t_{2n-1}$ and $t_{2n}:=t_0$).}
In short, we have approximations as
\begin{equation}\label{eq:def:p}
  V_k[\varphi]({\bf x}(t_m))\approx \sum_{\ell=-p}^{p}\alpha_{\ell,n'}^{\pm}({\bf x}(t_m))J_\ell(k|{\bf x}(t)-{\bf x}^{\pm}(t)|)e^{-i\ell\theta^{-}}
\end{equation}

%In order to achieve uniform errors using the quadratures~\eqref{eq:trap} for all indices $\ell$ (the integrands in equations~\eqref{eq:trap} get increasingly oscillatory and the near-singularities of the Hankel functions more stringent as the indices $\ell$ get larger), Fourier interpolation is used to oversample the density $\varphi$ on a finer uniform mesh of size $2\beta n$ where $\beta\in\mathbb{Z},\ \beta>1$ and the coefficients $\alpha_\ell^\pm$ are evaluated using the trapezoidal rule on the finer mesh.

We will also accurately evaluate the QBX expansion coefficients $\alpha_\ell^\pm$ using more general quadrature rules that are applicable in the case when $\Gamma$ is piecewise smooth. Specifically, we consider a panel representation of the boundary curve in the form $\Gamma=\bigcup_{{r}=1}^M \Gamma_{r}$ where the panels $\Gamma_{r}$ are non overlapping. In the case when $\Gamma$ exhibits corner points $\x_1,\x_2,\ldots,\x_P$, each corner point  is an end point of a panel. We thus have
\begin{equation}\label{eq:chebyshev:0}
  \alpha_\ell^\pm{(\x)}=\sum_{{r=1}}^M\alpha_{\ell,{r}}^\pm{(\x)},\qquad \alpha_{\ell,{r}}^\pm{(\x)}:=\frac{i}{4}\int_{\Gamma_{r}} H_\ell^{(1)}(k|\x^\pm-{\y}|)e^{i\ell(\theta')^\pm}\varphi({\y})\,{\rm d}{\y}.
\end{equation}
Assuming that each panel $\Gamma_{r}$ is parametrized in the form $\Gamma_{r}=\{{{\bf z}_r}(t): t\in[-1,1]\}$ where ${{\bf z}_r}:[-1,1]\to{\Gamma_r} $ is smooth,  we consider a Chebyshev mesh on the parameter space $[-1,1]$
\begin{equation*}
  {t_m}:=\cos \vartheta_m ,\quad \vartheta_m:=\frac{(2m-1)\pi}{2n_r},\quad m=1,\ldots,n_r.
\end{equation*}
The coefficients $\alpha_{\ell,r}^\pm$ are evaluated by a combination of oversampling and Fej\'er-Clenshaw-Curtis quadratures. Specifically,
\begin{equation}\label{eq:Clensh-Curtis}
{\alpha_{\ell,r}^\pm({\bf x}(t_m)}\approx{\alpha_{\ell,r,n'}^\pm({\bf x}(t_m))}:=\frac{i}{4}\sum_{j=1}^{n'_m}\omega_jH_\ell^{(1)}(k|
  {{\bf x}^\pm(t_m)}-{\bf z}_r(  {\tau_j})|)e^{i\ell(\theta')^\pm_j}\varphi({{\bf z}_r(\tau_j)})| {{\bf z}_r( \tau_j)}|
\end{equation}
with, as before, $n'_m =\beta n_m$ and consequently
\begin{equation*}
  {\tau_\ell}:=\cos \vartheta_\ell ,\quad \vartheta_\ell:=\frac{(2m-1)\pi}{2n'_r},\quad m=1,\ldots,n'_r.
\end{equation*}
The Fej\'er quadrature weights $\omega_j$, in turn, are given by
\begin{equation*}
  \omega_j:=\frac{2}{n_r'}\left(1-2\sum_{q=1}^{[n_r'/2]}\frac{1}{4q^2-1}\cos(2q\vartheta_j)\right),\quad j=1,\ldots,n'_r.
\end{equation*}
Again here, for a Chebyshev mesh point ${\bf z}_r(t_m)$ on $\Gamma$, we choose its corresponding centers along the exterior/interior normal to $\Gamma$ at ${\bf z}_r(t_m)$ located distance $\varepsilon_{j,r}=\min(|{\bf z}_r(t_m)-{\bf z}_r(t_{m+1})|,|{\bf z}_r(t_m)-{\bf z}_r(t_{m-1})|)$ from ${\bf z}(t_m)$. %(we assume that $\z(t_{m}$ and $\x_{j+1}$ are the grid points on $\Gamma$ adjacent to $\x_j$).
The evaluation of the density at $n_r'$ quadrature points is also carried out by (Chebyshev) interpolation of $\varphi$ from the coarse mesh with $n_r$ nodes.

% Again here, for a Chebyshev mesh point $\x_j$ on $\Gamma$ we chose its corresponding centers along the exterior/interior normal to $\Gamma$ at $\x_j$ located distance $\varepsilon(\x_j)=\min(|\x_j-\x_{j-1}|,|\x_j-\x_{j+1}|)$ from $\x_j$ (we assume that $\x_{j-1}$ and $\x_{j+1}$ are the grid points on $\Gamma$ adjacent to $\x_j$). We also use oversampling in the quadrature formulas~\eqref{eq:Clensh-Curtis},  that is, we apply those formulas on a finer Chebyshev mesh on $\Gamma_m$ of size $N_m=\beta n_m,\ \beta>1$ using Chebyshev interpolation to obtain the values of the density $\varphi$ from the coarse mesh with $n_m$ nodes onto the finer Chebyshev grid with $N_m$ nodes.

\section{Numerical results}\label{NR}

We present in this section a variety of numerical results about the accuracy of Nystr\"om discretizations of the elastodynamic BIE solvers discussed in this text. Specifically, we show far field accuracy results of solvers based on CFIE formulations as well as Helmholtz decomposition formulations~\eqref{eq:NH} and~\eqref{eq:DHBIE}. In addition, we study the iterative behavior of solvers based on the aforementioned formulations using GMRES~\cite{SaadSchultz} iterative solvers for the solution of the linear systems ensuing from Nystr\"om discretizations. While the size of the linear systems we considered allows for application of direct solvers, the iterative behavior of BIE formulations does shed light on the iterative properties of their three dimensional counterparts. Finally, we present numerical results concerning BIE based CQ solutions of time dependent elasticity scattering problems.

For a scattered elastic field ${\bf u}$ the associated longitudinal wave ${\bf u}_p$ and the transversal wave ${\bf u}_s$~ \cite{KupGeBa:1979} or \cite[Ch. 2]{AmKaLe:2009}  are defined as in equations~\eqref{eq:Hdecomp1}. The Kupradze radiation conditions~\cite{KupGeBa:1979} simply state that the functions $\varphi_p$ and $\varphi_s$ defined in equations~\eqref{eq:Hdecomp2} are radiative solutions of the Helmholtz equation in the unbounded domain $\Omega_+$ with {wave-numbers} $k_p$ and respectively $k_s$, which in vector form amounts to
\begin{equation}\label{eq:kupradze}
  {\bf u}_p(\x)=\frac{e^{ik_p|\x|}}{\sqrt{|\x|}}\left({\bf u}_{p,\infty}(\hat{{\x}})+\mathcal{O}\left(\frac{1}{|\x|}\right)\right)\qquad
  {\bf u}_s(\x)=\frac{e^{ik_s|\x|}}{\sqrt{|\x|}}\left({\bf u}_{s,\infty}(\hat{{\x}})+\mathcal{O}\left(\frac{1}{|\x|}\right)\right)
  \end{equation}
  as $|\x|\to\infty$ where $\hat{{\x}}=\x/|\x|$. We asses the accuracy of our solvers using the metric of maximum far field {errors} $\varepsilon_\infty$ of the quantities ${\bf u}_{p,\infty}$ and ${\bf u}_{s,\infty}$ evaluated at fine enough meshes on the unit circle $|\hat{{\x}}|=1$. Per usual, we consider both types of incident fields in our scattering experiments, that is, elastodynamic point sources and plane waves. The former type of incident field is used in the context of manufactured solutions wherein the numerical errors are evaluated against an exact solution. On the other hand, in the latter case of incident fields we computed far field errors with respect to reference solutions produced through very fine discretizations of the underlying BIE. We also report the number of unknowns $N$ used in the discretization of each of the two unknowns of in the elastodynamic BIE systems considered in this text (in the case of Nystr\"om discretizations that use the trigonometric polynomial space $\mathbb{T}_n$ we have $N=2n$). We start with numerical examples related to the method of manufactured solutions.

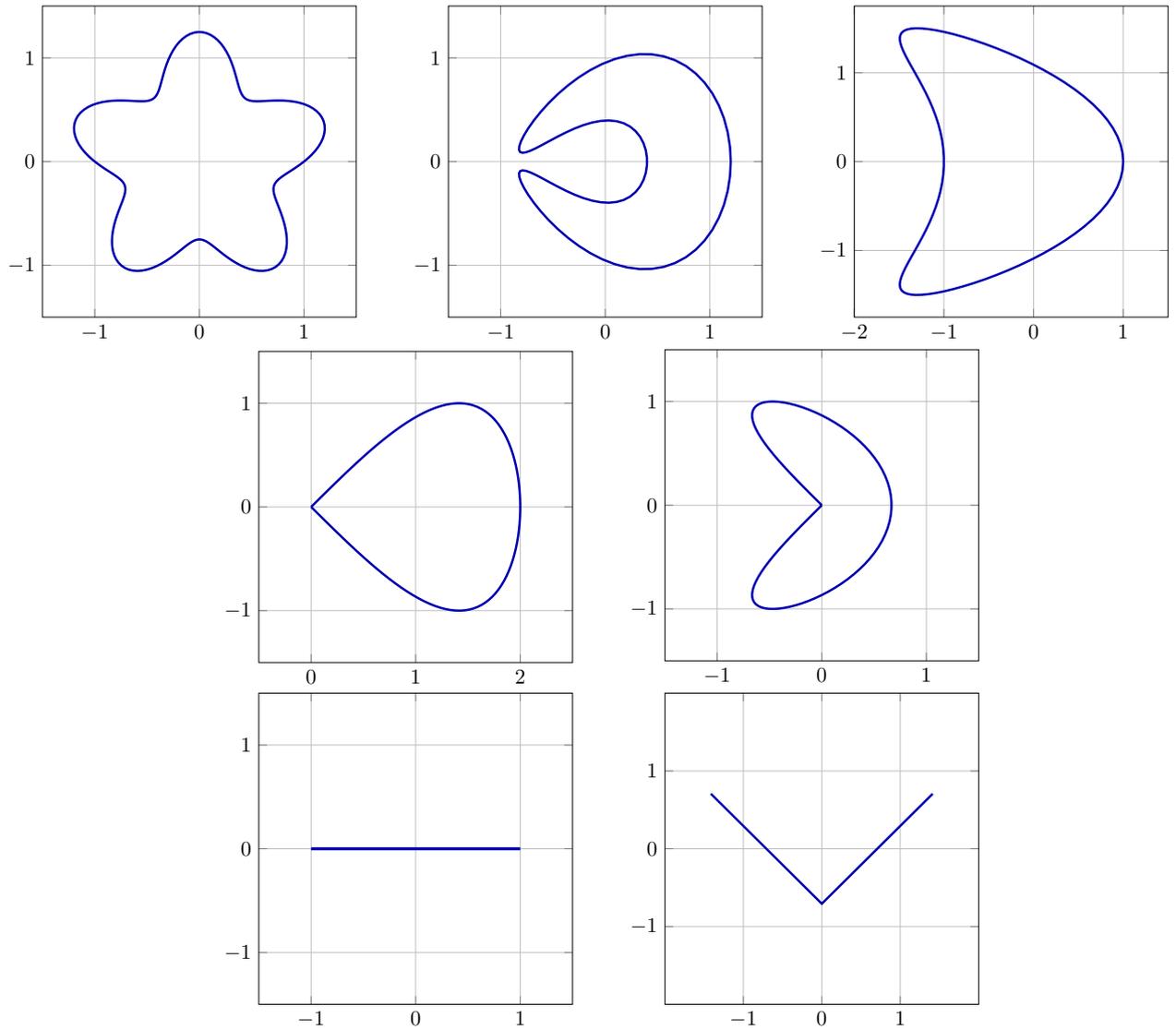
\begin{figure}
\begin{center}
%  \resizebox{0.31\textwidth}{!}{
% \begin{tikzpicture}
% [declare function = {f1(\x)=1); }]
% \begin{axis}[
%   xmin=-1.5 ,xmax=1.5 ,
%   ymin=-1.5,ymax=1.5,
%   grid=both,  axis equal image,
% xtick=    {-2, -1,0, 1,2 },
% ytick=    {-2, -1,0, 1,2 },
% ]
%   \addplot[domain=0:360,samples=200,variable=\x,very thick, blue!70!black,data cs=polar](x,{f1(x)});%
%    % {2.5 * (sin(deg(-5*t)))^2 * 2^(cos(deg(cos(deg(4.28*2.3*t)))))},%
%  %  t,t);
% \end{axis}
% \end{tikzpicture}
% }
%\quad
\resizebox{0.31\textwidth}{!}{
\begin{tikzpicture}
[declare function = {f1(\x)=1+ 0.25*sin(5*(\x) ); }]
\begin{axis}[
  xmin=-1.5 ,xmax=1.5 ,
  ymin=-1.5,ymax=1.5,
  grid=both,  axis equal image,
xtick=    {-2, -1,0, 1,2 },
ytick=    {-2, -1,0, 1,2 },
]
  \addplot[domain=0:360,samples=200,variable=\x,very thick, blue!70!black,data cs=polar](x,{f1(x)});%
   % {2.5 * (sin(deg(-5*t)))^2 * 2^(cos(deg(cos(deg(4.28*2.3*t)))))},%
 %  t,t);
\end{axis}
\end{tikzpicture}
}
\quad
 \resizebox{0.31\textwidth}{!}{
\begin{tikzpicture}
\begin{axis}[
  xmin=-1.5 ,xmax=1.5 ,
  ymin=-1.5,ymax=1.5,
  grid=both,  axis equal image,
xtick=    {-2, -1,0, 1,2 },
ytick=    {-2, -1,0, 1,2 },
]
  \addplot[domain=0:360,samples=100,variable=\t,very thick, blue!70!black](%
   % {2.5 * (sin(deg(-5*t)))^2 * 2^(cos(deg(cos(deg(4.28*2.3*t)))))},%
    {(cos( t)+2*cos(2* t))/2.5},
    {(sin( t)+ sin(2*t)+1/2* sin(3*t))/2-(-4* sin(t)+7*sin(2*t)-6*sin(3*t)+2*sin(4*t))/48}%
  );
\end{axis}
\end{tikzpicture}
}
\quad
 \resizebox{0.31\textwidth}{!}{
\begin{tikzpicture}
\begin{axis}[
  xmin=-2 ,xmax=1.5 ,
  ymin=-1.75,ymax=1.75,
  grid=both,  axis equal image,
xtick=    {-2, -1,0, 1,2 },
ytick=    {-2, -1,0, 1,2 },
]
  \addplot[domain=0:360,samples=100,variable=\t,very thick, blue!70!black](%
   % {2.5 * (sin(deg(-5*t)))^2 * 2^(cos(deg(cos(deg(4.28*2.3*t)))))},%
    {(cos( t)+0.65*cos(2* t)-0.65)},
    {1.5* sin(t)}%
  );
\end{axis}
\end{tikzpicture}
}
\\
\quad \resizebox{0.31\textwidth}{!}{
\begin{tikzpicture}
\begin{axis}[
  xmin=-0.5 ,xmax=2.5 ,
  ymin=-1.5,ymax=1.5,
  grid=both,  axis equal image,
xtick=    {0,1,2},
ytick=    {-2, -1,0, 1,2 },
]
  \addplot[domain=0:360,samples=100,variable=\t,very thick, blue!70!black](%
   % {2.5 * (sin(deg(-5*t)))^2 * 2^(cos(deg(cos(deg(4.28*2.3*t)))))},%
    {2*abs(sin(t/2))}, %
    {-sin(t)}%
  );
\end{axis}
\end{tikzpicture}
}
\quad \resizebox{0.31\textwidth}{!}{
\begin{tikzpicture}
\begin{axis}[
  xmin=-1.5 ,xmax=1.5 ,
  ymin=-1.5,ymax=1.5,
  grid=both,  axis equal image,
xtick=    {-1,0,1},
ytick=    {-2, -1,0, 1,2 },
]
  \addplot[domain=0:360,samples=100,variable=\t,very thick, blue!70!black](%
   % {2.5 * (sin(deg(-5*t)))^2 * 2^(cos(deg(cos(deg(4.28*2.3*t)))))},%
    {-2/3*(sin(3*t/2))}, %
    {-sin(t)}%
  );
\end{axis}
\end{tikzpicture}
}\\
\quad \resizebox{0.31\textwidth}{!}{
\begin{tikzpicture}
\begin{axis}[
  xmin=-1.5 ,xmax=1.5 ,
  ymin=-1.5,ymax=1.5,
  grid=both,  axis equal image,
xtick=    {-1, 0,1},
ytick=    {-1, 0,1 },
]
  \addplot[domain=-1:1,samples=5,variable=\t,ultra thick, blue!70!black](%
   % {2.5 * (sin(deg(-5*t)))^2 * 2^(cos(deg(cos(deg(4.28*2.3*t)))))},%
    {t}, %
    {0}%
  );
\end{axis}
\end{tikzpicture}
}
\quad \resizebox{0.31\textwidth}{!}{
\begin{tikzpicture}
\begin{axis}[
  xmin=-2 ,xmax=2 ,
  ymin=-2,ymax=2,
  grid=both,  axis equal image,
xtick=    {-1, 0,1},
ytick=    {-1, 0,1 },
]
  \addplot[domain=-1.414:1.414,samples=5,variable=\t,very thick, blue!70!black](%
   % {2.5 * (sin(deg(-5*t)))^2 * 2^(cos(deg(cos(deg(4.28*2.3*t)))))},%
    {t}, %
    {abs(t)-0.707}%
  );
\end{axis}
\end{tikzpicture}
}
\end{center}
\caption{\label{fig:geom} Geometries for the experiments considered in this section,  on the top, from left to right: the smooth domains, the  starfish,  the cavity, and the boomerang curve;  on the middle, the {\em corner} (Lipchitz) closed geometries:  the teardrop and the boomerang  curve; on the bottom, the open arcs:  the flat line and the V-shaped curve.}
\end{figure}

The geometries considered in our numerical experiments are:
\begin{enumerate}
 \item The {\em starfish} domain~\cite{hao2014high} whose $2\pi-$periodic paramaterization is given by
    \begin{equation}\label{hao2014high}
    {\bf x}(t)=\left(1+\frac{1}{4}\sin{5t}\right)(\cos{t},\sin{t}).
    \end{equation}

\item The {\em cavity}-like geometry whose parametrization is given by
\begin{multline}\label{eq:cavity}
{\bf x}(t) = (\tfrac14(\cos t+2\cos2t,A(t)/2-A_s(t)/48), \\
A(t)=\sin t+\sin 2t +1/2\sin 3t , A_s(t)=-4\sin t +7\sin 2t -6\sin 3t +2\sin 4t.                                                                                                                                                                                                                                                       \end{multline}

\item The {\em kite} domain cf. \cite{ColtonKress}
\[
{\bf x}(t)=(\cos t+0.65 (\cos 2 t-1),1.5 \sin t),
\]

\item The {\em teardrop} domain given by
\begin{equation}\label{eq:teardrop}
{\bf x}(t)=\left(2|\sin\tfrac{t}{2}|,-\sin{t}\right),
\end{equation}
which presents a corner point at ${\bf x}(0)=(0,0)$

\item The {\em boomerang} domain,
\[
{\bf x}(t)=  \left(-\tfrac{2}{3}\sin{\tfrac{3t}{2}}, -\sin{t}\right),
\]
which present also a corner point at  ${\bf x}(0)=(0,0)$.

\item The {\em flat} line, an open arc given by
\[
{\bf x}(t)=  (t,0),  \quad t\in[-1,1].
\]

\item The {\em V-shaped} polygonal line,
\[
{\bf x}(t)=  \left(t,|t|-\tfrac{1}{2\sqrt{2}}\right),  \quad t\in[-\tfrac{1}{\sqrt{2}},\tfrac{1}{\sqrt{2}}]
\]
and open curve with a corner point in $t=0$ of length 2.

\end{enumerate}
We depict in Figure \ref{fig:geom} such geometries.

  \subsection{The method of manufactured solutions}

  The method of manufactured solutions amounts to solving time-harmonic Navier equation with boundary value data produced by point sources $\x_0$ placed inside the scatter $ \Omega_-$
  \[
    {\bf u}^{\rm {inc}}(\x)=\Phi(\x,\x_0),\ \x\in\Gamma,\ \x_0\in \Omega_-
    \]
so that the solution of impenetrable scattering problems in the exterior domain $\Omega_+$ are the point sources themselves, that is, ${\bf u}(\x)=\Phi(\x,\x_0)$ for all $\x\in \Omega_+$. In the case of Dirichlet boundary conditions, looking for scattered fields in the form of a single and double layer potentials corresponding to boundary functional densities ${{\bm{\varphi}}}$ and respectively ${{\bm{g}}}$,  we solve the ensuing BIEs
    \[
    ({\bm V}{{\bm{\varphi}}})(\x)=\Phi(\x,\x_0)\quad{\rm and}\quad \frac{1}{2}{{\bm{g}}}(\x)+({\bm K}{{\bm{g}}})(\x)=\Phi(\x,\x_0),\ \x\in\Gamma.
    \]
 In the case of Neumann boundary conditions, the same approach leads to solving the following BIEs
   \[
    -\frac{1}{2}{{\bm{\varphi}}}(\x)+({\bm K}^\top{{\bm{\varphi}}})(\x)=T_{\x}\Phi(\x,\x_0), \quad{\rm and}\quad ({\bm W}{{\bm{g}}})(\x)=T_{\x}\Phi(\x,\x_0),\ \x\in\Gamma.
    \]
 In either type of boundary conditions we compare the numerical solutions in the far field against the exact point source solution ${\bf u}(\x)=\Phi(\x,\x_0)$. 
 The method of manufactured solutions in the case of Helmholtz decomposition BIE~\eqref{eq:DHBIE} and~\eqref{eq:NH} also amounts to consider the same incident field ${\bf u}^{\rm {inc}}(\x)=\Phi(\x,\x_0),\ \x\in\Gamma$ whose decomposition~\eqref{eq:Hdecomp2} is straightforward to effect by simply separating the $k_p$ and $k_s$ contributions via the appropriate Hankel functions in the fundamental solution of the Navier equation.
 
 In Table~\ref{comp5} we report the far field errors in the method of manufactured solutions in the case of a smooth starfish boundary.
Specifically, we present errors corresponding to (1) the single layer formulation with Dirichlet boundary conditions (in the rubric ${\bm V}$, given that the single layer BIO ${\bm V}$ is used to validate the method of manufactured solutions), (2) the double layer formulation with Dirichlet boundary conditions (in the rubric ${\bm K}$), and (3) the double layer formulation with Neumann boundary conditions (in the rubric ${\bm W}$). We mention that the BIO ${\bm K}^\top$ is the (real) $L^2\times L^2$ adjoint of the BIO ${\bm K}$, and at the discrete level the Nystr\"om matrices corresponding to those two operators are the transpose of one another. We considered two types of Nystr\"om discretizations, one based on Kussmaul-Martensen (K-M) logarithmic splitting strategy of all of the weakly singular kernels (which we {labeled} under the KM header) and another based on Alpert 10-th order quadrature (which does not require any splittings whatsoever of same kernels). The weights and the location of the off-grid nodes required by the Alpert 10-th order quadrature are tabulated in~\cite{alpert1999hybrid,hao2014high}, which, for the sake of brevity, we chose not to reproduce here. While the errors produced using the Alpert quadratures seem to saturate around $10^{-10}$ for both the single and double layer formulations of the Dirichlet problems and to $10^{-6}$ in the case of the hyper singular operators (in contrast, the corresponding errors using Kussmaul-Martensen splittings can reach full double precision levels), the application of Alpert quadratures to the elastodynamic 2D BIE solvers is significantly simpler.

\begin{table}
   \begin{center}
\begin{tabular}{|c|c|c|c|c|c|c|c|}
\hline
$\omega$ & $N$ & \multicolumn{2}{c|} {${\bm V}$} & \multicolumn{2}{c|} {${\bm K}$} & \multicolumn{2}{c|} {${\bm W}$} \\
\hline
 & & KM $\varepsilon_\infty$ & 10 Alpert $\varepsilon_\infty$ & KM $\varepsilon_\infty$ & 10 Alpert $\varepsilon_\infty$& KM $\varepsilon_\infty$ & 10 Alpert $\varepsilon_\infty$\\
\hline
16 & 64  & 3.0 $\times$ $10^{-2}$ & 1.5 $\times$ $10^{-2}$ &9.9$\times$ $10^{-2}$ &1.3 $\times$ $10^{-1}$ & 1.3 $\times$ $10^{-1}$&1.7$\times$ $10^{-1}$\\
16 & 128 & 7.1 $\times$ $10^{-7}$ & 1.7 $\times$ $10^{-6}$ & 1.1 $\times$ $10^{-3}$&3.1 $\times$ $10^{-5}$ & 1.6 $\times$ $10^{-6}$&7.4$\times$ $10^{-3}$ \\
16 & 256 & 5.5 $\times$ $10^{-14}$ & 1.7 $\times$ $10^{-10}$ & 4.2 $\times$ $10^{-13}$&3.4 $\times$ $10^{-9}$ & 1.6 $\times$ $10^{-9}$&3.2$\times$ $10^{-5}$\\
\hline
\hline
32 & 128  & 3.0 $\times$ $10^{-2}$ & 1.3 $\times$ $10^{-2}$ &2.5 $\times$ $10^{-1}$ &2.4 $\times$ $10^{-1}$ &3.4 $\times$ $10^{-1}$ &6.5$\times$ $10^{-2}$\\
32 & 256 & 2.2 $\times$ $10^{-7}$ & 1.4 $\times$ $10^{-6}$ & 2.4 $\times$ $10^{-4}$&7.0 $\times$ $10^{-6}$ &5.0 $\times$ $10^{-8}$ &4.6$\times$ $10^{-3}$\\
32 & 512 & 1.0 $\times$ $10^{-15}$ & 1.5 $\times$ $10^{-10}$ & 8.3 $\times$ $10^{-13}$&9.0 $\times$ $10^{-10}$ & 5.4 $\times 10^{-12}$ &2.8$\times$ $10^{-5}$ \\
\hline
\end{tabular}
\caption{Errors in the method of manufactured solution using the elastodynamics BIOs for the smooth starfish geometry for different values of the frequency $\omega$ and parameter values $\lambda=1$, $\mu=1$, at various levels of discretization using the Kussmaul-Martensen logarithmic splitting Nystr\"om discretizations as well as the 10th order Alpert discretization with parameters $a=6$ and $m=10$.\label{comp5}}
\end{center}
\end{table}

We present in Table~\ref{comp5qbx} errors in the method of manufactured solutions achieved by QBX Nystr\"om discretizations of the Helmholtz decomposition BIE formulations~\eqref{eq:DHBIE} of the elastodynamics scattering problem with Dirichlet boundary conditions (K-M results are available in the literature~\cite{dong2021highly} for the solution of the same equation~\eqref{eq:DHBIE}).  Both, equispaced and Chebyshev meshes are used in this experiment, with $p=16$, i.e. 33 terms in the expansion (see \eqref{eq:def:p}), $\beta=10$ in the oversampling and a global panel,  {i.e. $M=1$  which means that only one panel is being used in the Chebyshev mesh case}. We point out that in this experiment the one sided~\eqref{eq:qbxKtop} and two-sided~\eqref{eq:qbxKtopm} QBX expansions led to almost identical levels of accuracy for the same size of discretizations.

We continue in Table~\ref{comp5qbxPanel} with the same setup in the method of manufactured solutions but using different numbers $(M)$ of Chebyshev panels on the starfish contour (the interval $[0,2\pi]$ was split into $M$ equal parts and each panel corresponds to the mapping of one such subinterval via the parametrization of $\Gamma$) and values of the QBX expansion parameters $p$ and $\beta$. {The density is therefore computed at $n_r$ points per panel, and so it amounts to $n_r M$ on the curve. In this and the following experiments in this section, we will take $n_r=n$, i.e., for Chebyshev mesh the same number of points per panel will always be used. }

Qualitatively similar results are obtained when the double layer formulation~\eqref{eq:DDHBIE} (and hence the CFIE~\eqref{eq:DHCFIE}) is used for the solution of the Helmholtz decomposition approach.
\begin{table}[h]
   \begin{center}
\begin{tabular}{|c|c|c|c|}
\hline
$\omega$ & $N$ & \multicolumn{2}{c|} {Dirichlet b.c. Helmholtz decomposition BIE~\eqref{eq:DHBIE}} \\
\hline
 & & $\varepsilon_\infty$ Trapezoidal QBX~\eqref{eq:trap} & $\varepsilon_\infty$ Chebyshev QBX~\eqref{eq:Clensh-Curtis} \\
\hline
16 & 64  & 2.7 $\times$ $10^{-4}$ & 1.5 $\times$ $10^{-2}$ \\
16 & 128 & 6.6 $\times$ $10^{-8}$ & 1.7 $\times$ $10^{-6}$ \\
16 & 256 & 8.6 $\times$ $10^{-12}$ & 1.7 $\times$ $10^{-10}$\\ 
\hline
\hline
32 & 128  & 2.2 $\times$ $10^{-1}$ & 1.1 $\times$ $10^{-2}$ \\
32 & 256 & 5.3 $\times$ $10^{-5}$ & 1.8 $\times$ $10^{-7}$ \\
32 & 512 & 2.9 $\times$ $10^{-11}$ & 1.3 $\times$ $10^{-9}$ \\
\hline
\end{tabular}
\caption{Errors in the method of manufactured solution using the Helmholtz decomposition BIE~\eqref{eq:DHBIE} for the smooth starfish geometry with Dirichlet boundary conditions for different values of the frequency $\omega$ and parameter values $\lambda=1$, $\mu=1$, at various levels of discretization using QBX discretizations with global equispaced and Chebyshev meshes  with $M=1$ (only one panel) and with expansions parameters $p=16$ and $\beta=10$.\label{comp5qbx}}
\end{center}
\end{table}

\begin{table}[h]
   \begin{center}
\begin{tabular}{|c|c|c|c|c|c|}
\hline
\multicolumn{3}{|c|} {Chebyshev QBX~\eqref{eq:Clensh-Curtis}} &  \multicolumn{3}{c|} {Chebyshev QBX~\eqref{eq:Clensh-Curtis} } \\
\hline
 ${(M,{n})}$ & $(p,\beta)$ & $\varepsilon_\infty$ & ${(M,{n})}$ & $(p,\beta)$ & $\varepsilon_\infty$\\
 \hline
(8,16) & (4,4) & 2.4 $\times$ $10^{-3}$ & (32,8) & (2,2) & 1.5 $\times$ $10^{-2}$\\
(8,32) & (4,4) & 9.7 $\times$ $10^{-5}$ & (32,8) & (4,2) & 2.0 $\times$ $10^{-3}$\\
(8,32) & (8,6) & 4.0 $\times$ $10^{-6}$ & (32,8) & (6,4) & 2.6 $\times$ $10^{-4}$\\
(8,32) & (16,10) & 1.9 $\times$ $10^{-7}$ & (32,8) & (12,8) & 2.1 $\times$ $10^{-5}$\\
\hline
\end{tabular}
\caption{Errors in the method of manufactured solution using the Helmholtz decomposition BIE~\eqref{eq:DHBIE}  for the smooth starfish geometry with Dirichlet boundary conditions for the frequency $\omega=16$ and parameter values $\lambda=1$, $\mu=1$, using different numbers $(M)$ of Chebyshev panels of different sizes to discretize the boundary, as well as various choices of the QBX expansion parameters $p$ and $\beta$.
{The number of unknowns for each experiment is $N={n} M$}. \label{comp5qbxPanel}}
\end{center}
\end{table}

%We present in

Table~\ref{comp5qbxNeu} {illustrates} the accuracy levels achieved by the Nystr\"om K-M and QBX discretizations of the Neumann Helmholtz decomposition BIE~\eqref{eq:NH} in the case of the starfish scatterer.  The K-M Nystr\"om discretization is applied to the recasting of the Hessian operators that feature in equations~\eqref{eq:NH} via the Maue integration by parts techniques in formulas~\eqref{eq:jumpH} and~\eqref{eq:jumpHDL}. In the case of QBX discretizations we used a global equispaced mesh and corresponding trapezoidal quadratures for the evaluation of the Fourier-Bessel coefficients~\eqref{eq:trap}.  We discretized the BIOs that enter the BIE formulation~\eqref{eq:NH} via one sided QBX representations~\eqref{eq:qbxH1} and~\eqref{eq:qbxH2} and respectively two sided QBX representations~\eqref{eq:qbxH1m} and~\eqref{eq:qbxH2m}. Interestingly, it appears that the one sided QBX representations lead to more accurate solutions, which is {relevant} given that (a) they entail half the computational cost of two-sided expansions, and, more importantly (b) the one sided expansions are oblivious of the complicated jump relations~\eqref{eq:traceH}. We mention that we observed a similar slight increase in accuracy when using one sided QBX representations for other geometries and incident fields. We continue in Table~\ref{comp5qbxNeumannPanel} with numerical experiments concerning QBX discretizations of the BIE formulations~\eqref{eq:NH} using one sided expansions, and representations of $\Gamma$ as a union of Chebyshev panels and various levels of discretizations and values of the QBX expansion parameters. Again, high-order convergence is observed.  Qualitatively similar results are obtained when the double layer approach is used for the solution of the Neumann Helmholtz decomposition approach.

\begin{table}
   \begin{center}
\begin{tabular}{|c|c|c|c|c|}
\hline
$\omega$ & $N$ & \multicolumn{3}{c|} {Neumann b.c. Helmholtz decomposition BIE~\eqref{eq:NH}} \\
\hline
 & & K-M & One sided QBX $\varepsilon_\infty$ & Two sided QBX $\varepsilon_\infty$ \\
\hline
16 & 64  &  5.1 $\times$ $10^{-3}$ & 2.9 $\times$ $10^{-3}$ & 3.3 $\times$ $10^{-3}$\\ 
16 & 128 & 1.3 $\times$ $10^{-5}$ & 2.1 $\times$ $10^{-7}$ & 1.4 $\times$ $10^{-6}$ \\
16 & 256 & 8.2 $\times$ $10^{-11}$ & 2.8 $\times$ $10^{-11}$ & 1.2 $\times$ $10^{-10}$ \\
\hline
\hline
32 & 128  & 2.4 $\times$ $10^{-3}$& 1.4 $\times$ $10^{-3}$ & 1.3 $\times$ $10^{-3}$ \\
32 & 256 & 2.1 $\times$ $10^{-10}$& 4.9 $\times$ $10^{-8}$ & 5.0 $\times$ $10^{-8}$ \\
32 &  512 & 9.4 $\times$ $10^{-13}$ & 5.6 $\times$ $10^{-11}$ & 3.1 $\times$ $10^{-10}$ \\
\hline
\end{tabular}
\caption{Errors in the method of Neumann manufactured solution for the smooth starfish geometry using the Helmholtz decomposition formulation~\eqref{eq:NH} for different values of the frequency $\omega$ and parameter values $\lambda=1$, $\mu=1$, at various levels of discretization using K-M and QBX discretizations with a global equispaced mesh and expansion parameters $p=16$ and $\beta=10$. {The number of unknowns for each experiment is $N={n} M$}.\label{comp5qbxNeu}}
\end{center}
\end{table}

\begin{table}
   \begin{center}
\begin{tabular}{|c|c|c|c|c|c|}
\hline
\multicolumn{3}{|c|} {Chebyshev QBX~\eqref{eq:Clensh-Curtis}} &  \multicolumn{3}{c|} {Chebyshev QBX~\eqref{eq:Clensh-Curtis} } \\
\hline
 ${(M,{n})}$ & $(p,\beta)$ & $\varepsilon_\infty$ & ${(M,{n})}$ & $(p,\beta)$ & $\varepsilon_\infty$\\
 \hline
(8,16) & (4,4) & 2.3 $\times$ $10^{-3}$ & (32,8) & (2,2) & 1.3 $\times$ $10^{-2}$\\
(8,32) & (4,4) & 3.3 $\times$ $10^{-4}$ & (32,8) & (4,4) & 3.7 $\times$ $10^{-4}$\\
(8,32) & (6,6) & 6.6 $\times$ $10^{-6}$ & (32,8) & (8,6) & 3.5 $\times$ $10^{-5}$\\
(8,32) & (12,8) & 2.4 $\times$ $10^{-7}$ & (32,16) & (8,8) & 2.0 $\times$ $10^{-6}$\\
\hline
\end{tabular}
\caption{Errors in the method of manufactured solution using the Helmholtz decomposition BIE~\eqref{eq:NH}  for the smooth starfish geometry with Neumann boundary conditions for the frequency $\omega=16$ and parameter values $\lambda=1$, $\mu=1$, using different numbers $(M)$ of Chebyshev panels of different sizes to discretize the boundary, as well as various choices of the QBX expansion parameters $p$ and $\beta$. {The number of unknowns for each experiment is $N={n} M$}.\label{comp5qbxNeumannPanel}}
\end{center}
\end{table}

We will examine now the accuracy in the method of manufactured solutions in the case of Lipschitz geometries in Table~\ref{comp5c}. Specifically, we used the teardrop geometry \eqref{eq:teardrop} in our numerical tests in conjunction with polynomially graded meshes of order $p=4$ in the sigmoid transform and Kussmaul-Martensen kernel splittings as well as Alpert 3rd order method (using Lagrange interpolation with stencils of size 6 to evaluate the values of the densities at off grid locations). We have found that amongst all possible Alpert quadratures, the use of 3rd order Alpert quadratures delivers in practice nearly optimal accuracy given the lower order of regularity of BIOs functional densities in the Lipschitz case. However, we do not observe high-order convergence when we applied Alpert quadratures to the Nystr\"om discretization of the hyper singular operator $W$. We note that similar levels of accuracy are attained in the case of plane wave incidence.

\begin{table}
  \begin{center}
\begin{tabular}{|c|c|c|c|c|c|c|}
\hline
$\omega$ & $n$ & \multicolumn{2}{c|} {${\bm V}^w$} &  \multicolumn{2}{c|} {${\bm K}$} & ${\bm W}$ \\
\hline
 & & KM $\varepsilon_\infty$ & 3 Alpert KM $\varepsilon_\infty$ &{KM} $\varepsilon_\infty$ &  3 Alpert $\varepsilon_\infty$ & KM $\varepsilon_\infty$ \\
\hline
16 & 32  & 8.1 $\times$ $10^{-2}$ & 3.8 $\times$ $10^{-3}$ & 1.6 $\times$ $10^{-1}$ &  1.0 $\times$ $10^{-2}$& 1.7 $\times$ $10^{-1}$\\
16 & 64 & 8.4 $\times$ $10^{-5}$ & 6.4 $\times$ $10^{-5}$ & 3.4 $\times$ $10^{-4}$ & 2.3 $\times$ $10^{-3}$& 2.2 $\times$ $10^{-3}$ \\
16 & 128 & 4.1 $\times$ $10^{-12}$ & 2.0 $\times$ $10^{-6}$ & 2.3 $\times$ $10^{-9}$ & 2.8 $\times$ $10^{-4}$& 4.2 $\times$ $10^{-4}$ \\
\hline
\hline
32 & 64  & 4.7 $\times$ $10^{-2}$ & 6.4 $\times$ $10^{-3}$ & 2.0 $\times$ $10^{-1}$ & 1.1 $\times$ $10^{-2}$& 1.3 $\times$ $10^{-1}$ \\
32 & 128 & 8.1 $\times$ $10^{-6}$ & 1.7 $\times$ $10^{-4}$ & 1.7 $\times$ $10^{-4}$ & 1.7 $\times$ $10^{-3}$& 1.4 $\times$ $10^{-3}$\\
32 & 256 & 2.6 $\times$ $10^{-13}$ & 6.4 $\times$ $10^{-6}$ & 2.7 $\times$ $10^{-10}$ & 2.7 $\times$ $10^{-4}$ & 3.5 $\times$ $10^{-4}$ \\
\hline
\end{tabular}
\caption{Errors in the method of manufactured solution for the teardrop geometry for different values of the frequency $\omega$ and parameter values $\lambda=1$, $\mu=1$, at various levels of discretization using graded sigmoid meshes with $p=3$, the Kussmaul-Martensen logarithmic splitting Nystr\"om discretizations as well as the 3rd order Alpert discretization with parameters $a=2$ and $m=3$.\label{comp5c}}
\end{center}
\end{table}

  Tables~\ref{compTDcornerqbx} and~\ref{compTNcornerqbx} illustrate, again for Lipschitz scatterers, the accuracy of Helmholtz decomposition BIE formulations~\eqref{eq:DHBIE} and~\eqref{eq:NH} for the Navier equations with Dirichlet and respectively Neumann boundary conditions. Global Chebyshev meshes in those experiments with  Clenshaw-Curtis quadratures~\eqref{eq:Clensh-Curtis} are applied for the evaluations of the Fourier-Bessel coefficients in the QBX method. Again here, the one sided expansions appear to lead to more accurate solutions. We perform the same experiments in Tables~\ref{comp5qbxTDirichletPanel} and~\ref{comp5qbxTNeumannPanel} but using this time dyadic refinement of Chebyshev panels around the corner, and various values of the QBX expansion parameters $p$ and $\beta$. Similar levels of accuracy are observed for scatterers with re-entrant corners as well as multiple corners. We remark that the QBX discretizations of the BIE formulations~\eqref{eq:DHBIE} and~\eqref{eq:NH} are rather straightforward to implement even for Lipschitz boundaries, at least when one sided expansions are used (since one needs not take into account more complicated jump properties across $\Gamma$), and significantly simpler than discretizations of BIE based on the Navier fundamental solution.
\begin{table}
   \begin{center}
\begin{tabular}{|c|c|c|c|}
\hline
$\omega$ & $N$ & \multicolumn{2}{c|} {Dirichlet b.c. Helmholtz decomposition BIE~\eqref{eq:DHBIE}} \\
\hline
 & & $\varepsilon_\infty$ Chebyshev QBX one sided~\eqref{eq:Clensh-Curtis} & $\varepsilon_\infty$ Chebyshev QBX two-sided~\eqref{eq:Clensh-Curtis} \\
\hline
16 & 64  & 4.2 $\times$ $10^{-6}$ & 1.6 $\times$ $10^{-4}$ \\
16 & 128 & 6.8 $\times$ $10^{-9}$ & 1.6 $\times$ $10^{-6}$ \\
16 & 256 & 8.8 $\times$ $10^{-10}$ & 6.5 $\times$ $10^{-9}$\\ 
\hline
\hline
32 & 128  & 1.0 $\times$ $10^{-4}$ & 1.6 $\times$ $10^{-4}$ \\
32 & 256 & 8.5 $\times$ $10^{-8}$ & 4.6 $\times$ $10^{-7}$ \\
32 & 512 & 1.1 $\times$ $10^{-10}$ & 4.3 $\times$ $10^{-9}$ \\
\hline
\end{tabular}
\caption{Errors in the method of manufactured solution using the Helmholtz decomposition BIE~\eqref{eq:DHBIE} for the teardrop geometry with Dirichlet boundary conditions for different values of the frequency $\omega$ and parameter values $\lambda=1$, $\mu=1$, at various levels of discretization using QBX discretizations with global {equispaced} and Chebyshev meshes, with expansions parameters $p=16$ and $\beta=10$.\label{compTDcornerqbx}}
\end{center}
\end{table}

\begin{table}
   \begin{center}
\begin{tabular}{|c|c|c|c|}
\hline
$\omega$ & $N$ & \multicolumn{2}{c|} {Neumann b.c. Helmholtz decomposition BIE~\eqref{eq:NH}} \\
\hline
 & & $\varepsilon_\infty$ Chebyshev QBX one sided~\eqref{eq:Clensh-Curtis} & $\varepsilon_\infty$ Chebyshev QBX two-sided~\eqref{eq:Clensh-Curtis} \\
\hline
16 & 64  & 1.4 $\times$ $10^{-5}$ & 1.4 $\times$ $10^{-4}$ \\
16 & 128 & 3.6 $\times$ $10^{-6}$ & 2.0 $\times$ $10^{-5}$ \\
16 & 256 & 5.3 $\times$ $10^{-7}$ & 4.0 $\times$ $10^{-7}$\\ 
\hline
\hline
32 & 128  & 3.0 $\times$ $10^{-4}$ & 1.4 $\times$ $10^{-3}$ \\
32 & 256 & 2.0 $\times$ $10^{-6}$ & 4.0 $\times$ $10^{-7}$ \\
32 & 512 & 1.1 $\times$ $10^{-8}$ & 3.1 $\times$ $10^{-8}$ \\
\hline
\end{tabular}
\caption{Errors in the method of manufactured solution using the Helmholtz decomposition BIE~\eqref{eq:NH} for the teardrop geometry with Neumann boundary conditions for different values of the frequency $\omega$ and parameter values $\lambda=1$, $\mu=1$, at various levels of discretization using QBX discretizations with global {equispaced} and Chebyshev meshes, with expansions parameters $p=12$ and $\beta=8$.\label{compTNcornerqbx}}
\end{center}
\end{table}

\begin{table}
   \begin{center}
\begin{tabular}{|c|c|c|c|c|c|}
\hline
\multicolumn{3}{|c|} {Dirichlet BIE~\eqref{eq:DHBIE} QBX~\eqref{eq:Clensh-Curtis}} &  \multicolumn{3}{c|} {Dirichlet~\eqref{eq:DHBIE} QBX~\eqref{eq:Clensh-Curtis} } \\
\hline
 ${(M,{n})}$ & $(p,\beta)$ & $\varepsilon_\infty$ & ${(M,{n})}$ & $(p,\beta)$ & $\varepsilon_\infty$\\
 \hline
(8,16) & (4,4) & 7.8 $\times$ $10^{-4}$ & (32,16) & (2,2) & 1.6 $\times$ $10^{-2}$\\
(8,32) & (4,4) & 6.1 $\times$ $10^{-5}$ & (32,16) & (4,4) & 6.1 $\times$ $10^{-4}$\\
(8,32) & (6,6) & 8.1 $\times$ $10^{-6}$ & (32,16) & (8,6) & 5.0 $\times$ $10^{-5}$\\
(8,32) & (12,8) & 2.2 $\times$ $10^{-7}$ & (32,16) & (8,8) & 6.8 $\times$ $10^{-6}$\\
\hline
\end{tabular}
\caption{Errors in the method of manufactured solution using the Helmholtz decomposition BIE~\eqref{eq:DHBIE}  for the teardrop geometry with Dirichlet boundary conditions for the frequency $\omega=16$ and parameter values $\lambda=1$, $\mu=1$, using different numbers $(M)$ of Chebyshev panels with dyadic corner refinement of different sizes to discretize the boundary, as well as various choices of the QBX expansion parameters $p$ and $\beta$. {The number of unknowns for each experiment is $N={n} M$.}\label{comp5qbxTDirichletPanel}}
\end{center}
\end{table}

\begin{table}
   \begin{center}
\begin{tabular}{|c|c|c|c|c|c|}
\hline
\multicolumn{3}{|c|} {Neumann BIE~\eqref{eq:NH} QBX~\eqref{eq:Clensh-Curtis}} &  \multicolumn{3}{c|} {Neumann BIE~\eqref{eq:NH} QBX~\eqref{eq:Clensh-Curtis} } \\
\hline
 ${(M,{n})}$ & $(p,\beta)$ & $\varepsilon_\infty$ & ${(M,{n})}$ & $(p,\beta)$ & $\varepsilon_\infty$\\
 \hline
(8,16) & (4,4) & 1.9 $\times$ $10^{-3}$ & (32,16) & (2,2) & 1.0 $\times$ $10^{-2}$\\
(8,32) & (4,4) & 4.0 $\times$ $10^{-4}$ & (32,16) & (4,4) & 1.6 $\times$ $10^{-3}$\\
(8,32) & (6,6) & 7.7 $\times$ $10^{-6}$ & (32,16) & (8,6) & 3.7 $\times$ $10^{-4}$\\
(8,32) & (12,8) & 6.2 $\times$ $10^{-7}$ & (32,16) & (8,8) & 4.5 $\times$ $10^{-5}$\\
\hline
\end{tabular}
\caption{Errors in the method of manufactured solution using the Helmholtz decomposition BIE~\eqref{eq:NH}  for the teardrop geometry with Neumann boundary conditions for the frequency $\omega=16$ and parameter values $\lambda=1$, $\mu=1$, using different numbers $(M)$ of Chebyshev panels with dyadic corner refinement of different sizes to discretize the boundary, as well as various choices of the QBX expansion parameters $p$ and $\beta$. {The number of unknowns for each experiment is $N={n} M$.}\label{comp5qbxTNeumannPanel}}
\end{center}
\end{table}
In short, we have illustrated in this part the fact that QBX discretizations of the Helmholtz decomposition BIE formulations of the Navier scattering equations lead to the same levels of accuracy as the other Nystr\"om discretizations of the elastodynamics BIEs that use the Navier fundamental solution, while being simpler to implement. We turn our attention next to the iterative behavior of the BIE formulations considered in this text in the high frequency regime.

\subsection{Iterative behavior of BIE formulations}

We devote this section to comparisons between the iterative behavior of the various BIE formulations for the solution of high-frequency elastic impenetrable scattering problems considered in this text under plane wave incident fields of the form
\begin{equation}\label{eq:plwave}
  {\bf u}^{\rm inc}(\x)=\frac{1}{\mu}e^{ik_s\x\cdot{\bm d}}({\bm d}\times{\bm p})\times{\bm d}+\frac{1}{\lambda+2\mu}e^{ik_p\x\cdot{\bm d}}({\bm d}\cdot{\bm p}){\bm d}
\end{equation}
where the direction ${\bm d}$ has unit length $|{\bm d}|=1$. If the vector ${\bm p}$ is chosen such as ${\bm p}=\pm{\bm d}$, the incident plane is a pressure wave or P-wave. In the case when ${\bm p}$ is orthogonal to the direction of propagation ${\bm p}$, the incident plane wave is referred to as a shear wave or S-wave. We considered plane waves of direction ${\bm d}=\begin{bmatrix}0 & -1\end{bmatrix}^\top$ in all of our numerical experiments; in the case of S-wave incidence we selected ${\bm p}=\begin{bmatrix}1 &0\end{bmatrix}^\top$. We observed that other choices of the direction ${\bm d}$ and of the vector ${\bm p}$ lead to qualitatively similar results.

\subsubsection{Dirichlet boundary conditions}

We investigated the iterative behavior of Dirichlet integral solvers based on two formulations: (1) the CFIE formulation with the optimal coupling constant $\eta_D$ given in equation~\eqref{eq:eta_D}---which we refer to by the acronym ``CFIE $\eta_D$ opt''; and (2) the Helmholtz decomposition single layer BIE formulations~\eqref{eq:DHBIE} and their CFIE versions~\eqref{eq:DHCFIE}. We report in Tables~\ref{comp7}--~\ref{comp7b} the number of GMRES iterations required by each of these two BIE formulations to reach relative residuals of $10^{-5}$. The corresponding far field errors are also at the $10^{-5}$ level. in the case of smooth scatterers. Besides the starfish and the teardrop geometry, we also considered the Lipschitz boomerang geometry given by the discretization. As it can be seen from the results presented in Tables~\ref{comp7}-\ref{comp7b}, the Navier CFIE formulation~\eqref{eq:CFIER_D} with the optimal coupling parameter $\eta_D$ exhibits the best iterative behavior in the high-frequency regime. The iterations counts corresponding to the Navier CFIE formulation~\eqref{eq:CFIER_D} displayed in this section are almost identical for K-M and Alpert discretizations.

We note that given that the Navier double layer operators ${\bm K}$ and ${\bm K}^\top$ as well as the Helmholtz operators $\partial_sV_k$ and $W_k$ are not compact, neither of the Navier CFIE~\eqref{eq:CFIER_D} nor the single layer and CFIE Helmholtz decomposition BIE formulations~\eqref{eq:DHBIE} and~\eqref{eq:DHCFIE} are of the second kind. Nevertheless, the Navier CFIE formulation~\eqref{eq:CFIER_D} and the single layer Helmholtz decomposition BIE formulations~\eqref{eq:DHBIE} behave in practice as second kind formulations in the sense that the numbers of GMRES iterations required to reach a certain residual do not increase with more refined discretizations. However, that is not the case for the CFIE Helmholtz decomposition BIE formulations~\eqref{eq:DHCFIE}, largely on the account of the hyper singular operators $W_k$. Yet, the CFIE Helmholtz decomposition BIE formulations appear to be a superior alternative to the single layer Helmholtz decomposition BIE with regards to iterative solvers in the high frequency regime. We used QBX Nystr\"om discretizations based on global equispaced (in the smooth boundaries case) and respectively Chebyshev meshes for the discretizations of the Helmholtz decomposition BIE formulations~\eqref{eq:DHBIE}. We note that a Kussmaul-Martensen splitting Nystr\"om discretizations of the BIE formulations~\eqref{eq:DHBIE} is relatively straightforward to implement, and the ensuing numbers of GMRES iterations are slightly smaller than the ones resulting from QBX discretizations. In the QBX discretizations the use of one sided expansions led to identical results with respect to GMRES iterations with the versions using two sided expansions. Finally, we observed that the numbers of GMRES iterations grow with the number of panels used to represent the boundary $\Gamma$ when QBX discretizations are used (with the same size of overall discretization).

\begin{table}
   \begin{center} \renewcommand{\cellalign}{tc}
\begin{tabular}{|c|c|c|c|}
\hline
\makecell{$\omega$} & $N$ &\makecell{ \# iter CFIE~\eqref{eq:eta_D} $\eta_D$ opt}& \makecell{  \# iter Helmholtz decomposition\\
BIE~\eqref{eq:DHBIE} and~\eqref{eq:DHCFIE}}\\
\hline
\hline
10 & 64 & 25 & 205/100\\
20 & 128 & 27 & 302/124 \\
40 & 256 &  31 & 437/145 \\
80 & 512 &  34 &  727/166\\
160 & 1024 &  38 & 1401/183 \\
\hline
\end{tabular}
\caption{Numbers of GMRES iterations of various formulations to reach residuals of $10^{-5}$ for various BIE formulations of Dirichlet elastic scattering problems at high frequencies in the case when $\Gamma$ is the kite contour. The material parameters are $\lambda=2$ and $\mu=1$, and the incidence was P-wave. The discretizations used in these numerical experiments delivered results accurate at the level of $10^{-5}$.  The CFIE formulation was discretized using {Kussmaul-Martensen (K-M) logarithmic splittings Nystr\"om method}, while the BIE~\eqref{eq:DHBIE} were discretized using QBX with a global equispaced mesh and expansion parameters $p=12$ and $\beta=6$.\label{comp7}}
\end{center}
\end{table}

\begin{table}
   \begin{center} \renewcommand{\cellalign}{tc}
\begin{tabular}{|c|c|c|c|c|c|}
\hline
\makecell{$\omega$} & \makecell{$N$} &\makecell{\# iter CFIE~\eqref{eq:eta_D}\\
$\eta_D$ opt} & \makecell{\# iter one sided\\ QBX~\eqref{eq:DHBIE}} &\makecell{\# iter two sided\\
QBX~\eqref{eq:DHBIE}}&\makecell{\# iter one sided\\ QBX~\eqref{eq:DHCFIE}}\\
\hline
\hline
10  & 64   &  23 &   170 & 206    &  70  \\
20  & 128  &  28 &   292 & 338    &  80  \\
40  & 256  &  34 &   605 & 578    &  99  \\
80  & 512  &  40 &  1207 & 1137   & 137  \\
160 & 1024 &  50 &  1497 & 1405   & 186 \\
\hline
\end{tabular}
\caption{Numbers of GMRES iterations of various formulations to reach residuals of $10^{-5}$ for various BIE formulations of Dirichlet elastic scattering problems at high frequencies in the case when $\Gamma$ is the teardrop. The material parameters are $\lambda=2$ and $\mu=1$, and the incidence was P-wave. The discretizations used in these numerical experiments delivered results accurate at the level of $10^{-5}$. The BIE~\eqref{eq:DHBIE} and \eqref{eq:DHCFIE} were discretized using QBX with a global Chebyshev mesh and expansion parameters $p=12$ and $\beta=6$.
%{As comparison, the number of GMRES iterations required by the  method for CFIE formulations~\eqref{eq:DHCFIE} discretized by K-M Nystr\"om method} to reach the same GMRES residuals are 70, 80, 99,  137 and respectively 186.
\label{comp7t}}
\end{center}
\end{table}

\begin{table}
   \begin{center}             \renewcommand{\cellalign}{tc}
\begin{tabular}{|c|c|c|c|c|c|}
\hline
\makecell{$\omega$} & \makecell{$N$} & \makecell{\# iter CFIE~\eqref{eq:eta_D}\\ $\eta_D$ opt} & \makecell{\# iter one sided\\ QBX~\eqref{eq:DHBIE}} & \makecell{\# iter two sided\\ QBX~\eqref{eq:DHBIE}}&
\makecell{\# iter one sided\\ QBX~\eqref{eq:DHCFIE}}\\
\hline
\hline
10  &  128 &  29    &   198   &  204  &  72  \\
20  &  256 &  34    &   346   &  367  & 107\\
40  &  512 &  41    &   623   &  635  & 130\\
80  & 1024 &  49    &  1270   & 1229  & 165 \\
160 & 2048 &  59    &  1468   & 1418  & 214 \\
\hline
\end{tabular}
\caption{Numbers of GMRES iterations of various formulations to reach residuals of $10^{-5}$ for various BIE formulations of Dirichlet elastic scattering problems at high frequencies in the case when $\Gamma$ is the boomerang. The material parameters are $\lambda=2$ and $\mu=1$, and the incidence was P-wave. The discretizations used in these numerical experiments delivered results accurate at the level of $10^{-5}$. The BIE~\eqref{eq:DHBIE} and \eqref{eq:DHCFIE} were discretized using QBX with a global Chebyshev mesh and expansion parameters $p=8$ and $\beta=4$.
%The number of GMRES iterations required by the the Helmholtz decomposition CFIE formulations~\eqref{eq:DHCFIE} to reach the same GMRES residuals are 72, 107,  130,   165 and respectively 214.
\label{comp7b}}
\end{center}
\end{table}

We present next in Table~\ref{comp8a} results concerning scattering from arcs with Dirichlet boundary conditions. We considered a flat strip of length 2 and the V-shaped arc.
%of the same length where the line segments of the V-shaped geometry are orthogonal to each other.
We display  the number of iterations required by the classical single layer formulation (referred to as ${\bm V}^w$ in the column header) and the preconditioned on the left formulation that uses the operator composition ${\rm PS}_\kappa(Y_+){\bm V}^w$. We regularize the end-point square-root singularities of $T{\bf u}^{\rm tot}$~\cite{bruno2020regularized,bruno2019weighted} by resorting to its weighted version $T^w{\bf u}^{\rm tot}$ which, on account that the derivatives of the sigmoid transform vanish at the end-point of the $2\pi$ parametrization of the arc, vanishes itself at both endpoints $0$ and $2\pi$ and thus can be extended as a regular enough $2\pi$ periodic function. Clearly, the use of a weighted traction calls for the use of the weighted single layer BIO ${\bm V}^w$. We report numbers of GMRES iterations required to reach residuals of $10^{-6}$ and discretizations that deliver results at a $10^{-5}$ level of accuracy for the line segment and respectively $10^{-4}$ for the V-shaped arc. The ${\bm V}^w$ formulation is a first kind BIE, and hence the numbers of GMRES iterations required to reach a given tolerance grow upon refinement of discretizations, albeit modestly so. On the other hand, the analytically preconditioned formulation ${\rm PS}_\kappa(Y_+){\bm V}^w$ appears to behave like an integral equation of the second kind whose numbers of GMRES iterations are insensitive to the size of the dscretization for a given frequency. We remark that the computational overhead required by the application of the Fourier multiplier ${\rm PS}_\kappa(Y_+)$ is insignificant as its application can be performed efficiently using FFTs. The findings reported in Table~\ref{comp8a} appear to be qualitatively similar to those in~\cite{bruno2019weighted} where Calder\'on preconditioners are used.

\begin{table}
   \begin{center}  
\begin{tabular}{|c|c|c|c||c|c|}
\hline
$\omega$ & $N$ & \# iter ${\bm V}^w$& \# iter ${\rm PS}_\kappa(Y_+){\bm V}^w$& \# iter ${\bm V}^w$& \# iter ${\rm PS}_\kappa(Y_+){\bm V}^w$ \\
\hline
\hline
10 & 32/64 & 15/18 & 21/21 & 33/39 & 26/27\\
20 & 64/128 & 20/23 & 27/27 & 46/52 & 38/39\\
40 & 128/256 & 25/29 & 34/34 & 67/75 & 56/55\\
80 & 256/512 & 36/42 &  42/42 & 86/98 & 80/80\\
160 & 512/1024 & 51/61 & 53/53& 111/126 & 117/116\\
\hline
\end{tabular}
\caption{Numbers of GMRES iterations of various formulations to reach residuals of $10^{-6}$ for various BIE formulations of Dirichlet elastic scattering problems at high frequencies in the case when $\Gamma$ is a flat strip of length 2 (left panel) and a V-shaped arc of length 2 (right panel). The material parameters are $\lambda=2$ and $\mu=1$, and the incidence was an P-wave. The discretizations used in these numerical experiments delivered results accurate at the level of $10^{-5}$.\label{comp8a}}
\end{center}
\end{table}

\subsubsection{Neumann boundary conditions}

We present next in Tables~\ref{comp9}--\ref{comp12} results related to the iterative behavior of the BIE formulations considered in this paper for the solution of elastodynamic problems with Neumann boundary conditions in the case of smooth and Lipschitz scatterers. Specifically, we consider (1) the Navier CFIE formulation~\eqref{eq:CFIE_N} with the optimal coupling constant $\eta_N$, (2) the Navier CFIER formulation with the regularizing operator $\operatorname{R}^{\rm N}=({\rm PS}_\kappa(Y_+))^{-1}$ and (3) the single layer Helmholtz decomposition BIE formulation~\eqref{eq:NH} and its CFIE version. Both formulation~\eqref{eq:CFIE_N} and~\eqref{eq:NH} feature hyper singular operators and therefore the number of GMRES iterations required for the iterative solution of their associated Nystr\"om linear systems grows with the size of the discretizations.  The results reported in Tables~\ref{comp9}--\ref{comp12} correspond to K-M discretizations for smooth scatterers as well as one-sided QBX expansions for the discretization of the BIE~\eqref{eq:NH}; the use of two sided expansions leads to very similar iteration counts. Several remarks and observations are in order. First, the K-M and Alpert Nystr\"om discretizations of the Navier CFIE and CFIER formulations lead to almost identical iteration counts. Second, the Navier CFIE/CFIER formulations appear to enjoy superior iterative performance over the single layer and combined field Helmholtz decomposition formulations of the Neumann elastic scattering problems in the high frequency regime. Furthermore, the CFIE approach performs worse than the single layer approach in the case of the Helmholtz decomposition approach for Neumann elastic scattering, a fact which can be explained by the fact that the double layer formulation~\eqref{eq:jumpHDL} involves higher order derivatives than its single layer counterpart~\eqref{eq:traceH}. The construction of regularized formulations of the Helmholtz decomposition approach is currently under investigation.

\begin{table}
   \begin{center}  \renewcommand{\cellalign}{tc}
\begin{tabular}{|c|c|c|c|c|}
\hline
\makecell{$\omega$ }&\makecell{ $N$}  &\makecell{\# iter CFIE opt /\\ CFIER $\operatorname{R}^{\rm N}=({\rm PS}_\kappa(Y_+))^{-1}$ }& \makecell{\# Iter K-M /\\ one sided QBX~\eqref{eq:NH}} &\makecell{\# Iter one sided
QBX~\eqref{eq:NHCFIE}}
\\
\hline
\hline
10  &   128 &     51/29   &   187/260    &  238  \\
20  &   256 &     72/44   &   338/398    &  405 \\
40  &   512 &    121/71   &   620/737    &  783 \\
80  &  1024 &   195/123   &   1209/1332  & 1558 \\
160 &  2048 &   285/156   &   2085/2501  & 3459 \\
\hline
\end{tabular}
\caption{Numbers of GMRES iterations of various formulations to reach residuals of $10^{-5}$ for various BIE formulations of Neumann elastic scattering problems at high frequencies in the case when $\Gamma$ is the kite contour. The material parameters are $\lambda=2$ and $\mu=1$, and the incidence was an P-wave. The discretizations used in these numerical experiments delivered results accurate at the level of $10^{-4}$. Global equispaced meshes were used for QBX discretizations with expansion parameters $p=12$ and $\beta=6$.
%We mention that the iteration counts for the CFIE formulation of the Helmholtz decomposition approach using one sided QBX expansions are 238, 405, 783, 1558 and 3459.
\label{comp9}}
\end{center}
 \end{table}

\begin{table}
   \begin{center}           \renewcommand{\cellalign}{tc}
\begin{tabular}{|c|c|c|c|c|}
\hline
\makecell{$\omega$ }&\makecell{ $N$}  &\makecell{\# iter   CFIE opt / \\ CFIER $\operatorname{R}^{\rm N}=({\rm PS}_\kappa(Y_+))^{-1}$}& \makecell{\# Iter K-M / \\ one sided QBX~\eqref{eq:NH}} & \makecell{\# Iter one sided
QBX~\eqref{eq:NHCFIE}} \\
\hline
\hline
10  & 128   & 41/22   &    167/183  &  179 \\
20  & 256   & 52/31   &    220/237  &  279 \\
40  & 512   & 81/46   &    330/445  &  534 \\
80  &  1024 & 161/84  &    627/818  & 1061 \\
160 &  2048 & 348/167 &   1187/1419 & 2245 \\
\hline
\end{tabular}
\caption{Numbers of GMRES iterations of various formulations to reach residuals of $10^{-5}$ for various BIE formulations of Neumann elastic scattering problems at high frequencies in the case when $\Gamma$ is the starfish contour. The material parameters are $\lambda=2$ and $\mu=1$, and the incidence was an P-wave. The discretizations used in these numerical experiments delivered results accurate at the level of $10^{-4}$. Global equispaced meshes were used for QBX discretizations with expansion parameters $p=12$ and $\beta=6$.
%We mention that the iteration counts for the CFIE formulation of the Helmholtz decomposition approach using one sided QBX expansions are 179, 279, 534, 1061 and 2245.\label{comp9b}
}
\end{center}
 \end{table}

\begin{table}
   \begin{center}           \renewcommand{\cellalign}{tc}
\begin{tabular}{|c|c|c|c|c|}
\hline
\makecell{$\omega$ }&\makecell{ $N$} & \makecell{\# iter  CFIE opt /\\ DCFIER~\eqref{eq:S22}
$\operatorname{R}^{\rm N}=({\rm PS}_\kappa(Y_+))^{-1}$ }& \makecell{\# Iter one sided\\
QBX~\eqref{eq:NH}} &
\makecell{\# Iter one sided   \\
QBX~\eqref{eq:NHCFIE}}
\\
\hline
\hline
10  &  128  &    49/30  &  157  &  253\\
20  &  256  &    75/43  &  279  &  446\\
40  &  512  &   126/69  &  478  &  796\\
80  & 1024  &  211/116  &  910  & 1558 \\
160 & 2048  &  369/210  & 1876  & 1708\\
\hline
\end{tabular}
\caption{Numbers of GMRES iterations of various formulations to reach residuals of $10^{-5}$ for various BIE formulations of Neumann elastic scattering problems at high frequencies in the case when $\Gamma$ is a teardrop. The material parameters are $\lambda=2$ and $\mu=1$, and the incidence was an P-wave. The discretizations used in these numerical experiments delivered results accurate at the level of $10^{-4}$.  Global Chebyshev meshes were used for QBX discretizations with expansion parameters $p=12$ and $\beta=6$.
%We mention that the iteration counts for the CFIE formulation of the Helmholtz decomposition approach using one sided QBX expansions are 253, 446, 796, 1558 and 1708.
\label{comp11}}
\end{center}
\end{table}

\begin{table}
   \begin{center}           \renewcommand{\cellalign}{tc}
\begin{tabular}{|c|c|c|c|c|c|}
\hline
\makecell{$\omega$ }& \makecell{$N$} & \makecell{\# iter iter CFIE opt / \\ DCFIER~\eqref{eq:S22} $\operatorname{R}^{\rm N}=({\rm PS}_\kappa(Y_+))^{-1}$} & \makecell{\# Iter one sided\\ QBX~\eqref{eq:NH}} &
\makecell{\# Iter one sided   \\
QBX~\eqref{eq:NHCFIE}}
\\
\hline
\hline
10  &  128   &     72/36  &   216   &   252\\
20  &  256   &     97/50  &   319   &   435\\
40  &  512   &    137/77  &   650   &   798\\
80  & 1024   &   229/135  &  1191   &  1479 \\
160 & 2048   &   394/242  &  2267   &  2795 \\
\hline
\end{tabular}
\caption{Numbers of GMRES iterations of various formulations to reach residuals of $10^{-5}$ for various BIE formulations of Neumann elastic scattering problems at high frequencies in the case when $\Gamma$ is a boomerang. The material parameters are $\lambda=2$ and $\mu=1$, and the incidence was a P-wave.  The discretizations used in these numerical experiments delivered results accurate at the level of $10^{-4}$. Global Chebyshev meshes were used for QBX discretizations with expansion parameters $p=12$ and $\beta=6$.
%We mention that the iteration counts for the CFIE formulation of the Helmholtz decomposition approach using one sided QBX expansions are 252, 435, 798,  1479 and 2795.
\label{comp12}}
\end{center}
\end{table}

Finally, we conclude the numerical results with an illustration in Table~\ref{comp13} of scattering results in the case when $\Gamma$ is an open arc with Neumann boundary conditions. The solution ${\bf u}^{\rm tot}$ of the BIE~\eqref{eq:S12} can be shown to vanish at the end points of the (open) arc~\cite{chapko2000numerical,bruno2019weighted}, a salient feature that must be incorporated in the Nystr\"om discretization of the formulations~\eqref{eq:S12}. To do so, we assume the arc $\Gamma$ is smooth and follow the prescriptions in~\cite{chapko2000numerical}, that is, we use Chebyshev meshes on $\Gamma$ whereby the density function is parametrized on the interval $[0,\pi]$, and we take into account the fact that the density must vanish at the end points to extend it to the interval $[0,2\pi]$ via \emph{odd} extension. Just like in~\cite{chapko2000numerical}, we use the fact the operator ${\bm W}-{\bm W}_0$ is weakly singular as the basis of our Nystr\"om discretization of the hypersingular operator ${\bm W}$, and we project the discrete Nystr\"om equations back onto the interval $(0,\pi)$. We present results for the first kind formulation involving the hypersingular operator ${\bm W}$ (according to equation~\eqref{eq:S12}), as well as for the preconditioned formulation $({\rm PS}_\kappa(Y_+))^{-1}{\bm W}$. The discretization of the Fourier multiplier $({\rm PS}_\kappa(Y_+))^{-1}$ is performed using trigonometric interpolation for \emph{even} functions, that is, only cosines functions are used. The extension of the Helmholtz decomposition approach to the case of open arcs is currently underway.

\begin{table}
   \begin{center}  
\begin{tabular}{|c|c|c|c|}
\hline
$\omega$ & $n$ & \# iter ${\bm W}$~\eqref{eq:S12} & \# iter $({\rm PS}_\kappa(Y_+))^{-1}{\bm W}$\\
\hline
\hline
10 & 128 & 29 & 16\\
20 & 256 &  49 & 22\\
40 & 512 &  103 & 38\\
80 & 1024 & 225 & 59 \\
160 & 2048 & 490 & 76\\
\hline
\end{tabular}
\caption{Numbers of GMRES iterations of various formulations to reach residuals of $10^{-5}$ for various BIE formulations of Neumann elastic scattering problems at high frequencies in the case when $\Gamma$ is a flat strip of length 2. The material parameters are $\lambda=2$ and $\mu=1$, and the incidence was an P-wave. The discretizations used in these numerical experiments delivered results accurate at the level of $10^{-3}$.\label{comp13}}
\end{center}
\end{table}

\subsection{Time domain simulations}

We will show in this part how the numerical schemes introduced in this work  can be applied to solve the     transient elastic wave equation
 \begin{equation}\label{eq:td:01}
 \partial_t^2 {\bf u}(\x;t)=\diver \bm{\sigma}( {\bf u}(\x;t))\qquad{\rm in}\ {\Omega_+}\times(0,\infty)
 \end{equation}
 with either Dirichlet
\begin{subequations}\label{eq:td:02}
 \begin{equation}\label{eq:td:02a}
 {\bf u}(\x;0) ={\bf 0},\quad
{\bf u}(\x;t)={\bf f}(\x;t),\qquad \x\in\Gamma,\ t>0
 \end{equation}
 or Neumann conditions
  \begin{equation}\label{eq:td:02b}
 {\bf u}(\x;0) ={\bf 0},\quad
T{\bf u}(\x;t)={\bf g}(\x;t),\qquad \x\in\Gamma,\ t>0
 \end{equation}
 \end{subequations}
by   Convolution Quadrature (CQ)  Methods.
CQ {constitutes} a powerful framework for the numerical solutions of elastic wave equation~\cite{dominguez2015fully}. In particular, Runge-Kutta based CQ methods (RKCQ) can deliver higher order in time solvers for the wave equation~\cite{banjai2011runge} {than those based on multistep, as BDF, methods which are limited to order 2. We will present next the details of the implementation of  RKCQ  methods. The case of linear multi-step solvers (e.g. BDF2) was discussed in~\cite{dominguez2015fully}.}

\subsubsection{Runge Kutta Convolution Quadrature Methods for transient elastic wawe equation}

{We will follow closely the exposition in~\cite{betcke2017overresolving}  and we restrict oursevelves to the Dirichlet case \eqref{eq:td:02a}, since the Neumann problem cf. \eqref{eq:td:02b} is completely analogous}. First, problem \eqref{eq:td:01}-\eqref{eq:td:02a} is rewritten as
 as a first order system
 \begin{equation*}
    \begin{cases}
        \displaystyle \frac{\partial Y(\x;t)}{\partial t} = \mathcal LY(\x;t), \quad & (\x,t)\in (\mathbb{R}^2\setminus\Omega)\times(0,\infty)  \\
        BY(\x;t)=F(\x;t), \quad & (\x;t)\in\Gamma\times(0,\infty)\\
        Y(\x;0) = 0, \quad & \x \in \Omega_+
    \end{cases}
\end{equation*}
where we have introduced the following notations
\begin{align*}
    Y(\x;t) &=
    \begin{bmatrix}
    {\bf u}(\x;t) \\
    \displaystyle\frac{\partial {\bf u}(\x;t)}{\partial t } \end{bmatrix}
    &&\mathcal L = \begin{bmatrix}0 & I \\ \diver{\bm{\sigma}} & 0 \end{bmatrix} \\
    B &=\begin{bmatrix} I & 0 \\ 0 & 0 \end{bmatrix}   &&F(\x;t)=\begin{bmatrix}{\bf g}(\x;t)\\0 \end{bmatrix} .
\end{align*}
%CQ methods are time stepping methods that incorporate A-stable solvers for the solution of the first order system formulation of the wave equation and $\zeta$ transforms (basically discrete versions of Laplace transforms).
Given a final time $T$, the CQ methods deliver a sequence of $N+1$ approximations of the exact solution in $[0,T]$
%${\bf u}_{{\Delta t}}(\x;t_0), {\bf u}_{{\Delta t}}(\x,t_1),\ldots, {\bf u}_{{\Delta t}}(\x;T)$
\[
 [{\bf u}_{{\Delta t}}(\x; t_n)]_{0\le n\le N}, \quad t_n= n\Delta t.
\]
Here $N$ is a positive integer and $\Delta t= \frac{T}{N}$. We apply an $m$ stage Runge-Kutta (RK) scheme to the solution of the first order system reformulation of the elastic wave equation

\begin{equation}\label{eq:RK}
    \left\{\begin{array}{rcl}
        V_{i}(\x;t_n) &=& Y_{{\Delta t}}(\x;t_n)+ \Delta t\displaystyle \sum_{j=1}^m a_{ij}\mathcal LV_j(\x;t_n)\quad i\in{1,\ldots,m}, \\
        Y_{{\Delta t}}(\x;t_{n+1}) &=& Y_{{\Delta t}}(\x;t_n) + \Delta t \displaystyle\sum_{j=1}^m b_j\mathcal LV_j(\x;t_n)
    \end{array}\right.
\end{equation}
where the ensemble of matrices and vectors
\begin{equation*}
    {\bm{A}}= [a_{ij}]_{1\leq i,j\leq m }\quad
    {\bm{b}} = [b_j]_{1\leq j \leq m} \quad
    {\bm{c}}=[c_j]_{1\leq j\leq m},
\end{equation*}
make up the Butcher Tableau associated with the RK scheme. The A-stability requirement on the underlying ODE solver of CQ methods motivates the choice of \emph{stiffly accurate} {or $L-$stable} RK schemes~\eqref{eq:RK}. In these schemes,
\begin{equation}\label{eq:interstellar:01}
\bm{b}^\top \bm{A}^{-1} = \begin{bmatrix}
                 0&0&\cdots&0&1
                \end{bmatrix},
\end{equation}
{which trivially implies $Y_{{\Delta t}}(\x;t_{n+1})=V_{m}(\x;t_n)$. Hence we can focus on the computation of the internal stages of the RK method}.

Applying the $\zeta-$transform to the equations~\eqref{eq:RK}, {and taking into account \eqref{eq:interstellar:01}}, we obtain
\begin{equation}
    \left\{\begin{array}{rcl}
        V_{i}(\x;\zeta) &=& Y_{{\Delta t}}(\x;\zeta)+ \Delta t\displaystyle \sum_{j=1}^m a_{ij}\mathcal LV_j(\x;\zeta)\quad i\in{1,\ldots,m}

         \\
         \frac{Y_{{\Delta t}}(\x;\zeta)}{\zeta} &=&   Y_{{\Delta t}}(\x;\zeta) + \displaystyle \Delta t\sum_{j=1}^m  b_j\mathcal{L}V_j(\x;\zeta).
     \end{array}\right.
    \label{z_intermediate}
\end{equation}
where
\[
 V_{i}(\x;\zeta):=\sum_{n\geq 0}V_i(\x;t_n)\zeta^n,\quad
 Y_{d}(\x;\zeta):=\sum_{n\geq 0}Y_{{\Delta t}}(\x;t_n)\zeta^n.
\]
We note that the second equation in~\eqref{z_intermediate} leads to
\begin{equation}
    Y_{{\Delta t}}(\x;\zeta) = \frac{\zeta}{1-\zeta}\Delta t\sum_{j=1}^m b_j\mathcal LV_j(\x;\zeta).
\end{equation}
% , which, in turn, imposes the condition that the last row of the matrix $A$ must coincide with the vector $b$ {and that $c_m=1$}. That is,   Consequently, it holds that
% \begin{equation*}
%     V_m(\x;\zeta) = \frac{Y_{{\Delta t}}(\x;\zeta)}{\zeta},
% \end{equation*}
%which, in turn,

Writing the vector quantities $V_j$ in explicit form as $[R_j(\x;\zeta)\ S_j(\x;\zeta)]^\top$ and plugging it   in the first equation in~\eqref{z_intermediate},  yields to the following system of equations
\begin{equation}\label{eq:syst_1}
   \begin{bmatrix}
   R_i(\x;\zeta)\\
   S_i(\x;\zeta)
   \end{bmatrix}= \Delta t\sum_{j=1}^m \left ( \frac{\zeta}{1-\zeta}b_j + a_{ij} \right )
   \mathcal{L} \begin{bmatrix}
   R_j(\x;\zeta)\\
   S_j(\x;\zeta)
   \end{bmatrix} ,\quad i=1,\ldots,m.
\end{equation}
{Gathering the quantities above in vector form}
\[
  \bm{R}(\x;\zeta) := [R_i(\x;\zeta)]_{i=1,\ldots,m},\quad
  \bm{S}(\x;\zeta) := [S_i(\x;\zeta)]_{i=1,\ldots,m},
\]
we have {from the definition of the matrix differential operator ${\cal L}$}
\begin{align*}
    \bm{R}(\x;\zeta) &=
    \Delta t\left( \frac{\zeta}{1-\zeta}{\bf 1}\bm{b}^\top + \bm{A}\right)\otimes\bm{S}(\x;\zeta), &
    \bm{S}(\x;\zeta)& =
    \Delta t\left( \frac{\zeta}{1-\zeta}{\bf 1}\bm{b}^\top + \bm{A}\right)\otimes
  \diver \bm{\sigma}(\bm{R}(\x;\zeta))
\end{align*}
(the differential operator $\diver \bm{\sigma}$ is obviously applied element-wise) where
\[
 \mathbf{1} = [1,\ldots,1]^\top \in\mathbb{R}^m
\]
and $\bm{A}\otimes \bm{B}$ denoting the Kronecker product:
\[
 \bm{A}\otimes \bm{B} =
 \begin{bmatrix}
  a_{11} \bm B& \cdots &a_{1m}\bm B\\
  \vdots & \ddots&\vdots\\
  a_{m1} \bm B& \cdots &a_{mm}\bm B
 \end{bmatrix}   .
\]

% \begin{align*}
%     R_i(\x;\zeta) &= \Delta t\sum_{j=1}^m \left ( \frac{\zeta}{1-\zeta}b_j + a_{ij} \right ) S_j(\x;\zeta) \\
%     S_j(\x;\zeta) & =\Delta t\sum_{\ell=1}^m \left ( \frac{\zeta}{1-\zeta}b_j + a_{j\ell} \right )  \diver{\bm{\sigma}}R_\ell(\x;\zeta).
% \end{align*}
%Let us denote $\mathcal R(\x;\zeta) := [R_1(\x;\zeta),\ldots,R_m(\x;\zeta)]$. The last two equations above can then be written in compact form as a system of Helmholtz equations
{Therefore, we can recast the linear system~\eqref{eq:syst_1} in the equivalent form}
\begin{equation}
    \left( \frac{\Delta(\zeta)}{\Delta t}  \right )^2{\otimes \bm{R}}(\x;\zeta) =  \mathcal \diver{\bm{\sigma}}({\bm{R}}(\x;\zeta))
    \label{coupled}
\end{equation}
where $\Delta(\zeta)$ is a matrix operator defined by
\begin{equation*}
    \Delta(\zeta) = \left( {\bm{A}}+\frac{\zeta}{1-\zeta}\mathbf{1}  \bm{b}^{{\top}}     \right )^{-1}.
\end{equation*}
A natural idea is to attempt to decouple the linear system~\eqref{coupled} via diagonalization. To this end, we begin by diagonalizing the operator
$\Delta(\zeta)$. Let ${\bm{P}(\zeta)}=[P_{ij}(\zeta)]_{1\le i,j\le m}$ be the matrix consisting of the eigenvectors of $\Delta(\zeta)$ and define
\[
{
\bm{D}(\zeta)={\rm diag}(\gamma_1(\zeta),\ldots,\gamma_m(\zeta))}
\]
the diagonal matrix with the corresponding eigenvalues of $\Delta(\zeta)$. Thus
\[
 \Delta(\zeta)={\bm{P}(\zeta)}{\bm{D}(\zeta)}{\bm{P}(\zeta)}^{-1}.
\]
Clearly, the system~\eqref{coupled} is decoupled into $m$ frequency domain Navier equations
\begin{equation}\label{eq:decoupled}
    \left(\frac{\gamma_j(\zeta)}{\Delta t} \right)^2 W_j(\x;\zeta) = \diver{\bm{\sigma}} (W_j(\x;\zeta)),\ 1\leq j\leq m
\end{equation}
where
\[
W_j(\x;\zeta)=\sum_{\ell=1}^m  ({P}^{-1}(\zeta))_{j\ell}(\x;\zeta)R_\ell(\x;\zeta),\quad\text{or equivalently}\quad \bm{R}(\x;\zeta) =
\bm{P}(\zeta)\otimes  \bm{W}(\x;\zeta)  .
\]
Finally, boundary conditions ought to be provided for the frequency domain Navier equations~\eqref{eq:decoupled}.
%First we provide boundary conditions for the intermediate RK stages $V_j$; these take on the form
{For the intermediate RK stages $V_j$ these take on the form}
\begin{align*}
    BV_j &= F(\x;t_n+c_j\Delta t) = \begin{bmatrix} {\bf g}(\x;t_n+c_j\Delta t)\\0 \end{bmatrix} ,\quad
    {\bf g}(\x;t  + \bm{c}\, \Delta t) =
    \begin{bmatrix}
     g(\x;t  + c_1\, \Delta t)\\
     g(\x;t  + c_2\, \Delta t)\\
     \vdots\\
     g(\x;t  + c_m\, \Delta t)
    \end{bmatrix}.
\end{align*}
Applying the $\zeta-$transform to the boundary conditions above we get
\begin{equation}
  \gamma_\Gamma   \bm{R} (\x;\zeta) =\bm{G}(\x;\zeta)
 {=}
    \sum_{n \geq 0}{\bf g}(\x;t_n + \bm{c}\, \Delta t)\zeta^n
    %\quad \x\in \Gamma
\end{equation}
and thus we derive the corresponding boundary conditions for the Navier solutions $W_j$:
\begin{equation}
   % W_j(\x;\zeta) = \sum_{\ell=1}^m{ (P^{-1}(\zeta))_{j\ell}}G_\ell(\x;\zeta),\quad \x\in\Gamma.
   \gamma_\Gamma \bm{W}(\x;\zeta) =\bm{w}(\x;\zeta)= \begin{bmatrix}
                                                     w_1(\x;\zeta)\\
                                                     w_2(\x;\zeta)\\
                                                     \vdots\\
                                                     w_m(\x;\zeta)
                                                    \end{bmatrix}
:=
\bm{P}^{-1}(\zeta)\otimes  \bm{G}(\x;\zeta).
\end{equation}
Consequently, the following modified Navier equations must be solved in the Laplace domain
\begin{equation}\label{eq:CQh2}
    \begin{cases}
        \diver{\bm{\sigma}}(W_j(\,\cdot\,;\zeta)) - \left ( \frac{\gamma_j(\zeta)}{\Delta t} \right)^2 W_j(\,\cdot\,;\zeta) =0 \quad &\text{in } \Omega_+ \\
        \gamma_\Gamma W_j(\,\cdot\,;\zeta)=
        w_j(\,\cdot\,;\zeta).
        %\sum_{\ell=1}^m{ (P^{-1}(\zeta))_{j\ell}}G_\ell(\x;\zeta)\quad \x\in\Gamma.
    \end{cases}
\end{equation}
Once the Navier equations above are solved, the $\zeta$-transform {approximation of} the solution of the wave equation is retrieved via the formula (recall that $Y_{{\Delta t}}(\x;t_{n+1})=V_{m}(\x;t_n)$)
\begin{equation}
   \bm{u}_{\Delta t}(\x;\zeta) = \zeta R_m(\x;\zeta) = \zeta\sum_{j=1}^m {P}_{mj}(\zeta)W_j(\x;\zeta).
\end{equation}
Finally, reverting to the physical domain from the Laplace domain is performed in the same manner as in the case of linear multistep methods:
\[
  \bm{u}_{\Delta t}(\x;n \Delta t) := \frac{R^{-n}}{N+1} \sum_{\ell = 0}^N  \bm{u}_{\Delta t}(\x; R^\ell\zeta_{N+1}^\ell ) \zeta_{N+1}^{\ell n}\approx \bm{u}(\x,n \Delta t),\qquad \zeta_{N+1} =\exp\left(\tfrac{2\pi i}{N} \right).
\]
Here $0<R<1$ with $R=\varepsilon^{\tfrac{1}{2N+2}}$ with $\varepsilon$ being the unit round-off as suggested optimal value.

% We will use in our experiments two-stage and three-stage Implicit RK Radau IIa methods. The two-stage RK Radau II a method gives rise to third order in time CQ and will be denoted by the acronym IRK3. The three-stage RK Radau II a method gives rise to fifth order in time CQ and will be denoted by the acronym IRK5. We note that the higher order in time convergence that can be achieved by RK CQ  methods entail commensurately more solutions of frequency domain modified Navier equations. We refer to the ensemble of Navier equations~\eqref{eq:CQh2} as modified Navier equations because they entail complex rather than real frequencies. We use BIE to solve the ensemble of Navier boundary value problems~\eqref{eq:CQh2} and thus we obtain a numerical solution of the elastic wave equation~\eqref{eq:td} via the CQ methodology. Nystr\"om discretizations of BIE formulations of the elastodynamics boundary value problems~\eqref{eq:CQh2} face a distinct challenge arising from the fact that the Bessel functions used in the splitting of the functions $\phi_j, j=0,1$ (and hence throughout the singularity splitting of all four BIOs of Navier equations) grow exponentially for {\em complex} arguments of large absolute values. While this situation can be mitigated by the use of cut-off functions~\cite{turc1} to localize the splitting procedure around singularities, the details become cumbersome in the Navier case, and we advocate for the use of Alpert quadratures that are oblivious to those issues.

\subsubsection{Numerical results for Convolution Quadrature methods}

In our experiments the boundary conditions~\eqref{eq:td:02b}  on $\Gamma\times(0,\infty)$ correspond to the incident field
\[
{\bf u}^{\rm inc}(\x;t)=H(c_Lt-\x\cdot\bm{d})\sin(c_Lt-\x\cdot\bm{d})\bm{d},\ c_L
\]
where $\bm{d}=(1,0)$, $\lambda=1,\ \mu=1$ and $H$ is a smoothed version of the Heaviside function. We present numerical results for two smooth scatterers, the starfish   \eqref{hao2014high} and the cavity-like~\eqref{eq:cavity} geometry.
both Dirichlet and Neumann boundary conditions, and we present numerical results in the near field (the observation points are placed equispaced on a circle situated at distance 1 from the scatterers) for a final time $T=3$. The numerical approximations ${\bf u}_{{\Delta t}}(\cdot\ ; t_n)$ of time domain simulations are produced at the time grid $t_n=n\Delta t,\ 0\leq n\leq N$ such that $T=N\Delta t=3$
with BDF2 (multistep) or RK3/RK5 (multistage) solvers.
The ensemble of Laplace domain Navier equations~\eqref{eq:CQh2} are solved using the single layer potential formulation (that is, we use the BIO $\bm{V} $) in the case of Dirichlet boundary conditions (e.g. Table~\ref{comp5TDc1}) and  the double layer formulation (that is, we use the BIO $\bm{W}$) in the case of Neumann boundary conditions (e.g. Table~\ref{comp5TDc2}). These BIE of the first kind are the ones that are most frequently used in the CQ literature~\cite{dominguez2015fully}, and therefore we chose to present numerical results based on those in order to illustrate the levels of accuracy than can be achieved by our Nystr\"om discretizations. The discretization of these first kind BIE, in turn, is effected through Nystr\"om discretizations based on Alpert quadratures of order 10 (i.e. $a=6$ and $m=10$) using $2n=512$ discretization points on the boundaries $\Gamma$ for each of the Laplace domain frequencies that feature in equations~\eqref{eq:CQh2}. We present in Table~\ref{comp5TDc1} and Table~\ref{comp5TDc2} the orders of convergence in time (under the headings estimated order of convergence "e.o.c.") achieved by the CQ BDF2 and CQ RK3 when the errors $\varepsilon_\infty$ are computed in the near field at final time $T=3$ with respect to reference solutions (in the near field)  obtained using CQ RK5 for $N=2048$ time steps. We note that the CQ errors in {this} case appear to saturate at the level of $10^{-7}/10^{-8}$ when RK3 solvers are used, which is the highest level of accuracy the CQ methods can actually achieve~\cite{banjai2010multistep}. We end this section with a comment on the iterative behavior of BIE formulations of the Laplace domain Navier problems~\eqref{eq:CQh2}. Based on our experience, it is the BIE formulations of the second kind that perform best in this regard: in the case of the starfish/cavity geometry, the double layer formulation in the case of Dirichlet boundary conditions and the single layer formulation in the case of Neumann boundary conditions require at most $29/32$ iterations to reach GMRES residuals of $10^{-7}$ for all the frequencies corresponding to all values of time steps $N$ considered, with similar levels of accuracy to those reported in Table~\ref{comp5TDc1} and Table~\ref{comp5TDc2}. In contrast, both formulations considered in Tables~\ref{comp5TDc1} and~\ref{comp5TDc2}, being first kind formulations, require larger numbers of iterations for convergence---up to one order of magnitude more, as the number of frequencies is increased, a feature that is, shared by CFIER formulations as well. Qualitatively similar results are obtained when the Helmholtz decomposition integral formulations~\eqref{eq:DHBIE} and~\eqref{eq:NH} are used to solve the ensemble of CQ Laplace domain problems.

\begin{table}
  \begin{center}
\begin{tabular}{|c|c|c|c|c|c|c|c|c|}
\hline
$N$ & \multicolumn{4}{c|} {Starfish} &  \multicolumn{4}{c|} {Cavity} \\
\hline
 & BDF2 $\varepsilon_\infty$ & e.o.c.  &RK3 $\varepsilon_\infty$ &  e.o.c. & BDF2 $\varepsilon_\infty$ & e.o.c.  & RK3 $\varepsilon_\infty$ &  e.o.c.\\
\hline
32 & 5.5 $\times$ $10^{-2}$ &  & 1.1 $\times$ $10^{-3}$ &  & 2.4 $\times$ $10^{-2}$ &   & 1.1 $\times$ $10^{-3}$ &  \\
64 & 1.8 $\times$ $10^{-2}$ & 1.57 & 1.4 $\times$ $10^{-4}$ & 2.97 & 7.5 $\times$ $10^{-3}$ & 1.73 & 1.3 $\times$ $10^{-4}$ & 3.01\\
128 & 5.0 $\times$ $10^{-3}$ & 1.90 & 1.7 $\times$ $10^{-5}$ & 3.00 & 1.9 $\times$ $10^{-3}$ & 1.98 & 1.7 $\times$ $10^{-5}$ & 3.00\\
256 & 1.3 $\times$ $10^{-3}$ & 1.94 & 2.1 $\times$ $10^{-6}$ & 3.00  & 4.9 $\times$ $10^{-4}$ & 1.93 & 2.2 $\times$ $10^{-6}$ & 3.00\\
512 & 3.2 $\times$ $10^{-4}$ & 2.02 & 2.7 $\times$ $10^{-7}$ & 3.01 &1.3 $\times$ $10^{-4}$ & 1.92 & 2.8 $\times$ $10^{-7}$ & 2.97 \\
1024 & 8.0 $\times$ $10^{-5}$ & 1.99 & 3.3 $\times$ $10^{-8}$ & 3.01 & 4.1 $\times$ $10^{-5}$ & 1.65 & 2.7 $\times$ $10^{-7}$ & 3.01\\
2048 & 2.0 $\times$ $10^{-5}$ & 2.00 & 2.3 $\times$ $10^{-8}$ & &1.2 $\times$ $10^{-5}$ & 1.77 & 2.9 $\times$ $10^{-7}$& \\
\hline
\end{tabular}
\caption{CQ simulations of the solution of equations~\eqref{eq:td:01} with Dirichlet boundary conditions \eqref{eq:td:02a} using Nystr\"om discretizations based on Alpert 10th order quadratures for the single layer BIE formulations that use the BIO $V$ of the ensemble of Laplace domain Navier equations~\eqref{eq:CQh2}.\label{comp5TDc1}}
\end{center}
\end{table}

\begin{table}
  \begin{center}
\begin{tabular}{|c|c|c|c|c|c|c|c|c|}
\hline
$N$ & \multicolumn{4}{c|} {Starfish} &  \multicolumn{4}{c|} {Cavity} \\
\hline
 & BDF2 $\varepsilon_\infty$ & e.o.c.  &RK3 $\varepsilon_\infty$ &  e.o.c. & BDF2 $\varepsilon_\infty$ & e.o.c.  & RK3 $\varepsilon_\infty$ &  e.o.c.\\
\hline
32 & 1.9 $\times$ $10^{-2}$ &  & 1.4 $\times$ $10^{-3}$ & & 2.6 $\times$ $10^{-3}$ &  & 3.9 $\times$ $10^{-4}$ &  \\
64 & 6.5 $\times$ $10^{-3}$ & 1.58 & 1.7 $\times$ $10^{-4}$ & 3.01 & 1.1 $\times$ $10^{-3}$ &  1.24 & 4.8 $\times$ $10^{-5}$ & 3.01\\
128 & 1.7 $\times$ $10^{-3}$ &  1.93 & 2.2 $\times$ $10^{-4}$ & 2.86&3.6 $\times$ $10^{-4}$ &  1.61 & 5.9 $\times$ $10^{-6}$ & 3.03\\
256 & 4.0 $\times$ $10^{-4}$ & 2.06 & 2.1 $\times$ $10^{-5}$ & 3.07 &1.0 $\times$ $10^{-4}$ &  1.81 & 8.9 $\times$ $10^{-7}$ & 2.73\\
512 & 9.8 $\times$ $10^{-5}$ & 2.03 & 2.4 $\times$ $10^{-6}$ & 3.06 &2.7 $\times$ $10^{-5}$ &  1.92 & 2.6 $\times$ $10^{-7}$ & 1.77 \\
1024 & 2.4 $\times$ $10^{-5}$ & 2.01 & 2.9 $\times$ $10^{-7}$ & 3.08 &7.0 $\times$ $10^{-6}$ & 1.93  & 2.6 $\times$ $10^{-7}$ & \\
2048 & 6.1 $\times$ $10^{-6}$ & 2.00 & 3.2 $\times$ $10^{-7}$ & &1.9 $\times$ $10^{-6}$ &  1.87 & 2.5 $\times$ $10^{-7}$ & \\
\hline
\end{tabular}
\caption{CQ simulations of the solution of equations~\eqref{eq:td:01} with Neumann boundary conditions \eqref{eq:td:02b} using Nystr\"om discretizations based on Alpert 10th order quadratures for the double layer BIE formulations that use the BIO $W$ of the ensemble of Laplace domain Navier equations~\eqref{eq:CQh2}.\label{comp5TDc2}}
\end{center}
\end{table}

\section{Conclusions}
   
We presented two high-order Nystr\"om methods for the discretization of the four {BIO}s associated with time-harmonic Navier equations in two dimensions for both smooth as well as Lipschitz boundaries. These discretizations were used for the solution of elastic scattering problems in frequency and time domain based on BIE formulations.  We presented high-order Nystr\"om discretizations of the Helmholtz BIOs that feature in the Helmholtz decomposition BIE alternative formulation of elastodynamics scattering problems. Comparisons between the iterative behavior of various BIE formulations of elastodynamics scattering problems were carried out in the high-frequency regime. Extensions to three dimensional configurations are currently underway.

\section*{Acknowledgments}
 Catalin Turc gratefully acknowledges support from NSF through contract DMS-1908602.

\bibliography{biblioEL14062022}

\end{document}